\documentclass[twoside,leqno]{article}

\usepackage[letterpaper]{geometry}

\usepackage{siamproceedings}

\usepackage[T1]{fontenc}
\usepackage{amsfonts}
\usepackage{graphicx}
\usepackage{epstopdf}
\usepackage{enumitem}

\newenvironment{keywords}{\textbf{Keywords.}}{}
\usepackage{algorithmic}
\ifpdf
  \DeclareGraphicsExtensions{.eps,.pdf,.png,.jpg}
\else
  \DeclareGraphicsExtensions{.eps}
\fi

\newsiamremark{remark}{Remark}
\newsiamremark{hypothesis}{Hypothesis}
\crefname{hypothesis}{Hypothesis}{Hypotheses}
\newsiamthm{claim}{Claim}

\usepackage{amsopn}
\DeclareMathOperator{\diag}{diag}

\usepackage{xcolor}
\newcommand{\R}{\mathbb{R}}
\DeclareMathOperator{\dom}{dom}
\newcommand{\epsg}{\epsilon_g}
\newcommand{\epsH}{\epsilon_H}
\newcommand{\cF}{{\cal F}}
\newcommand{\cN}{{\cal N}}

\newcommand{\cS}{{\cal S}}
\newcommand{\cB}{{\cal B}}
\newcommand{\bfone}{{\mathbf{1}}}
\newcommand{\sigmin}{{\sigma_{\min}}}
\newcommand{\sigmax}{{\sigma_{\max}}}
\newcommand{\Dx}{\Delta x}
\newcommand{\Dl}{\Delta \lambda}
\newcommand{\Ds}{\Delta s}

\newcommand{\cK}{{\cal K}}

\begin{document}

\newcommand\relatedversion{}
\renewcommand\relatedversion{\thanks{The full version of the paper can be accessed at \protect\url{https://arxiv.org/abs/0000.00000}}} 

\title{\Large Optimization in Theory and Practice \relatedversion}
    \author{Stephen J. Wright\thanks{Department of Computer Sciences, University of Wisconsin-Madison, Madison, WI (\email{swright@cs.wisc.edu}, \url{https://wrightstephen.github.io/sw_proj/}).}}

\date{}

\maketitle


\fancyfoot[R]{\scriptsize{Copyright \textcopyright\ 20XX by SIAM\\
Unauthorized reproduction of this article is prohibited}}





\begin{abstract} 
Algorithms for continuous optimization problems have a rich history of design and innovation over the past several decades, in which  mathematical analysis of their convergence and complexity properties plays a central role.
Besides their theoretical properties, optimization algorithms are interesting also for their practical usefulness as computational tools for solving real-world problems.
There are often gaps between the practical performance of an algorithm and what can be proved about it.
These two facets of the field --- the theoretical and the practical --- interact in fascinating ways, each driving innovation in the other. 
This work focuses on the development of algorithms in two areas --- linear programming and unconstrained minimization of smooth functions --- outlining major algorithm classes in each area along with their theoretical properties and practical performance, and highlighting how advances in theory and practice have influenced each other in these areas. In discussing theory, we focus mainly on non-asymptotic complexity, which are upper bounds on the amount of computation required by a given algorithm to find an approximate solution of problems in a given class.
\end{abstract}


\section{Introduction.} \label{sec:intro}

Optimization is the process of minimizing a real-valued function of many variables, subject to restrictions or constraints on the values of those variables. 
This simply stated paradigm encompasses mathematical formulations for an enormous range of applications in many areas of science, engineering, and economics.
It also admits numerous refinements into classes and subclasses of problems (e.g. linear programming, nonsmooth optimization, integer programming), each of which provides formulations for many applications and each of which can be studied mathematically.
These classes have been refined over the years in response to the demands of application areas.
Over the past 15 years, for example, machine learning has been a rich source of optimization problems characterized by many variables, a great deal of data, and certain structures in the objective functions that can be exploited by algorithms.

Mathematics is used in optimization to study the properties of problem {\em formulations} and to devise {\em algorithms} and analyze their properties.
These strands are intertwined. 
The study of formulations includes the existence and characterization of solutions; the stability of solutions under perturbations of the problem data; and  duality, which provides an alternative perspective on a given formulation that can shed light on the problem and can sometimes be exploited by algorithms.
Algorithms for optimization problems are typically iterative, generating a sequence of candidate solutions that (ideally) converges to a solution for the problem at hand. 
Algorithms are rigorous in the sense that they allow certain claims about their properties to be proved, but within the rigorous framework, they contain heuristics to guide certain decisions, for example, the choice of pivot in a simplex algorithm for linear programming from among a (potentially enormous) list of valid choices.
Along with  mathematical expertise and imagination, a good deal of intuition and domain knowledge for the key application areas is required for successful algorithm design, adaptation, and application.
Development of {\em software} to implement these algorithms requires many other areas of expertise to be brought to bear, involving broader mathematical knowledge and software engineering expertise. 
The mathematical issues include the vagaries of finite-precision arithmetic, robustness in the face of ill-posed formulations, the interface with numerical linear algebra, and parallel implementation, among many others.
Engineering issues (some of which can also be analyzed mathematically) include management of hierarchical memory and data movement. 

The focus of this paper is the analysis of algorithms for two classes of {\em continuous} optimization problems --- problems in which the variables are vectors or matrices involving real numbers that are allowed to vary continuously, not being restricted to certain discrete values such as the integers.
We consider the convergence and complexity properties of these methods.
Are the iterates guaranteed to converge to points with the characteristics of a solution?
Can we bound the amount of computation needed to find an approximate solution over all problems in a class --- and how do we measure ``computation,'' given that the cost of obtaining information about functions can vary widely from one application to the next?
Can we {\em lower-bound} the amount of computation needed for algorithms of a certain type on the hardest problems in a certain class?

We focus especially on the interplay between the theoretical properties of algorithms and their practical performance, which  has driven a great deal of research in optimization over the past several decades.
For some problem classes, the practical performance of methods on typical problems in the class mirrors the theory well.
In other classes, there is a wide gap between theory and practice, with practical performance being much faster on typical problems than the theory predicts.
Such gaps can sometimes  be narrowed by refining the theory, perhaps by narrowing the problem class in an interesting way.
A recent development of this type has occurred in the minimization of smooth nonconvex functions (see \Cref{sec:nonconvex}).
In its full generality, this problem is intractable, but many subclasses of nonconvex problems have been identified (particularly in machine learning) for which global minima can be found efficiently. 
In other cases of wide gaps between theory and practice, further analysis and computation can improve our understanding of the nature of the gap. 
For example, we may ask: Are the loose bounds provided by the theory due to rare worst-case instances? Can we quantify the rarity of these instances? Can we mollify these instances, for example by showing that a nearby instance is usually easy to solve?
Gaps may prompt a search for new kinds of algorithms with better theoretical and possibly even better practical properties. 
A striking instance of this type has occurred in linear programming (see \Cref{sec:lp}), where the search for algorithms with polynomial complexity led eventually to the development of interior-point methods, which have proved to be important in practice for solving large-scale problems.

A paper on a similar theme to this one, but in the domain of computer science, was presented by Donald Knuth at the World Computer Congress in 1989 \cite{Knu91}.
It contains many enlightening observations about the relationship between theory and practice in that field, summarizing them in the quote ``The best theory is inspired by practice and the best practice is inspired by theory.'' 
Optimization has adhered to this dictum throughout its history as both a  mathematical and highly practical discipline.

\subsection*{Outline.}
We consider two classes of problems in this paper: (1) linear programming and (2) unconstrained optimization of a smooth real-valued function. 
There is a wide variety of algorithms for both classes, with rich mathematical properties and extensive computational experience.
\Cref{sec:opt} defines these problems and describes their optimality conditions, which are algebraic conditions that characterize solutions.
We discuss convergence and complexity analysis in broad terms in \Cref{sec:complexity},  outlining the reasons for the observed gaps between this theory and computational experience with these methods.
Linear programming is the focus of \Cref{sec:lp}, while unconstrained optimization is tackled in \Cref{sec:f}. \Cref{sec:discussion} contains some concluding remarks.

Space prevents us from considering other areas of optimization in which is interesting interplay between theory and practice, including constrained optimization, optimization of finite-sum functions (see \Cref{eq:finite-sum} below), and parallel methods.

\subsection*{Notation.}

We use $\R^n$ to denote the Euclidean space of dimension $n$.
For $x,y \in \R^n$, $\| x \|$ denotes the Euclidean norm and $\langle x,y\rangle$ denotes the inner product $x^Ty$. 
Components of a vector $x \in \R^n$ are denoted by subscripts: $x_i$, $i=1,2,\dotsc,n$.
The vector inequality $x \ge y$ for $x, y \in \R^n$ means that $x_i \ge y_i$, $i=1,2,\dotsc,n$.
When $A,B \in S\R^{n \times n}$ (symmetric $n \times n$ real matrices), we use $A \succeq B$ to indicate that $A-B$ is positive semidefinite. $A \succ B$ means that $A-B$ is positive definite.
$f:\R^n \to \R$ is a convex function when its epigraph is a convex set, equivalently, $f(\alpha x+(1-\alpha)y) \le \alpha f(x) + (1-\alpha) f(y)$ for all $x,y \in \dom f$ and $\alpha \in [0,1]$.

The order notation $O(\alpha)$ for some $\alpha>0$ means  that the quantity in question is bounded by some constant multiple of $\alpha$.  $\tilde{O}(\alpha)$ means that there might be an additional factor that depends logarithmically on such quantities as  dimension or solution accuracy.

The cardinality of  finite set $B$ is denoted by $|B|$.
We use $\bfone$ to denote a vector whose entries are all $1$.

\section{Formulations and Optimality Conditions.} \label{sec:opt}

Optimization formulations provide mathematical descriptions of practical problems around which theory and algorithms can be developed.
Any given practical problem can be formulated in many different ways, all ``valid'' in the sense that their solutions are equivalent and achieve the goals of the application.
However, some of these formulations may be much harder to solve than others, due to redundancies or degeneracies (and many other more subtle reasons). 
Expertise is needed to choose the formulation that can be solved most efficiently with the algorithmic tools at hand.
Conversely, designers of algorithms and software try to make their tools as robust as possible to such problematic aspects of formulations.

{\em Optimality conditions} are characterizations of solutions of a formulation, usually in terms of algebraic conditions in a form that can serve as ``target" for algorithms. 
Algorithms use these conditions to recognize when they are at or near a solution, and often use them to guide progress toward solutions as well. 

We describe several formulations in this section which will be studied in the remainder of the paper.

\subsection{Unconstrained Optimization.}
We consider first the following  unconstrained optimization problem:
\begin{equation} \label{eq:f}
\min_{x \in \R^n} \, f(x),
\end{equation}
where $\R^n$ denotes real vectors of length $n$ and $f:\R^n \to \R$ is a smooth real-valued function (at least Lipschitz continuously differentiable).
Finding (global) solutions of \cref{eq:f} is intractable for general smooth nonconvex $f$, although algorithms of branch-and-bound type can be applied when we have additional knowledge about $f$ (for example, a bound on its Lipschitz constant) and when we can restrict the search to a finite domain. 
Exceptions occur when $f$ is a convex function; the solutions of most convex problems can be found efficiently. 
In general, algorithms for \cref{eq:f} seek {\em local solutions}, which are points $x^*$ such that there is a radius $r>0$ such that $f(x^*) \le f(x)$ for all $x$ with $\|x-x^* \| \le r$.
(A {\em global solution} has $r=\infty$ in this definition.)
Local solutions can be characterized by {\em optimality conditions}. 
When $f$ is differentiable, the {\em first-order necessary} condition for $x^*$ to be a local solution of \cref{eq:f} is $\nabla f(x^*)=0$. 
(Points satisfying this condition are termed {\em stationary}.)  
When $f$ is a convex function, this condition is sufficient for $x^*$ to be a {\em global} solution. 
When $f$ is twice continuously differentiable, {\em second-order necessary conditions} for $x^*$ to be a local solution of \cref{eq:f} are that $\nabla f(x^*)=0$ and  $\nabla^2 f(x^*) \succeq 0$. {\em Second-order sufficient conditions} replace the second condition by positive definiteness: $\nabla^2 f(x^*) \succ 0$.

Since most algorithms that target stationary points or solutions converge only asymptotically to such points, we define {\em approximate} optimality conditions, allowing these algorithms to terminate finitely when such conditions are satisfied. 
An approximate first-order condition could be  $\|\nabla f(x^*) \| \le \epsg$, for some small $\epsg>0$, while an approximate second-order condition could add the condition  $\nabla^2 f(x^*) \succeq -\epsH I$, for some small $\epsH>0$.

Other notions of approximate optimality include the {\em optimality gap} $f(x)-f^*$, where $f^*$ is the value of $f$ at the solution, and $\mbox{\rm dist} (x,\cS)$, where $\cS$ is the set of solutions to \cref{eq:f}.
Neither condition is ``checkable" in general, as we need to know the solutions before evaluating them.
Nevertheless they are useful in the sense that it may be possible to obtain a bound on the computational effort required to satisfy these conditions approximately to within some positive tolerance $\epsilon$, that is, $f(x)-f^* \le \epsilon$ and $\mbox{\rm dist} (x,\cS) \le \epsilon$, respectively.


\subsection{Linear Programming.}
Linear programming (LP) is the problem class that jump-started the modern era of optimization in the 1940s. 
In LP, the objective function and the algebraic equalities and inequalities defining the feasible set are all affine.
Without loss of generality, LP can be written in the following {\em standard form}:
\begin{equation}
    \label{eq:lp}
    \min_{x \in \R^n} \, c^Tx \;\; \mbox{subject to} \;\; Ax=b, \;\; x \ge 0,
\end{equation}
where $A \in \R^{m \times n}$ is the constraint matrix, $b \in \R^m$ is the right-hand side of the constraints, and $c \in \R^n$ is the cost vector.
We assume throughout that $A$ has full row rank. (An elementary preprocessing step can ensure that this property holds.)
There are comprehensive theories  of optimality and duality for LP.
The problem \cref{eq:lp} may not have a solution --- it can be infeasible (when there is no vector $x$ satisfying the constraints $Ax=b$, $x \ge 0$) or unbounded (when there is a sequence of  vectors $x^k$, all satisfying the constraints, for which $\lim_{k \to \infty} \, c^Tx^k = -\infty$).
A vector $x$ is a solution to \cref{eq:lp} if and only if there is a vector pair $(\lambda,s) \in \R^m \times \R^n$ such that 
\begin{equation} \label{eq:lp.opt}
Ax=b, \quad A^T \lambda+s=c, \quad (x,s) \ge 0, \quad x_i s_i = 0, \; i=1,2,\dotsc,n.
\end{equation}
When such conditions are satisfied, $(\lambda,s)$ is a solution to another linear program, called the {\em dual} of \cref{eq:lp}:
\begin{equation}
    \label{eq:lp.dual}
    \max_{(\lambda,s) \in \R^m \times \R^n} \, b^T\lambda \;\; \mbox{subject to} \;\; A^T\lambda+s=c, \;\; s \ge 0.
\end{equation}
It is easy to show that for $(x,\lambda,s)$ satisfying \cref{eq:lp.opt}, we have $c^Tx-b^T\lambda = x^Ts = 0$, so the optimal values of the primal and dual LPs are the same, a property known as {\em strong duality}. 
When the first three conditions in \cref{eq:lp.opt} are satisfied, but not necessarily the final condition, we have that $x$ is a feasible point for \cref{eq:lp} and $(\lambda,s)$ is feasible for \cref{eq:lp.dual}, and $c^Tx \ge b^T\lambda$, a property known as {\em weak duality}. 
Note that the dual form \cref{eq:lp.dual} can be stated equivalently as
\begin{equation}
    \label{eq:lp.dual2}
    \max_{\lambda \in \R^m} \, b^T\lambda \;\; \mbox{subject to} \;\; A^T\lambda \le c.
\end{equation}

We can define approximate solutions of the linear programs \cref{eq:lp} and \cref{eq:lp.dual} in terms of approximate satisfaction of the optimality conditions \cref{eq:lp.opt}. 
One particular measure to which we refer in \Cref{sec:lp} asssumes that the feasibility conditions (the first three conditions in \cref{eq:lp.opt}) are satisfied exactly, while the quantity $\mu$ defined by
\begin{equation}
    \label{eq:mu}
    \mu := \frac{1}{n} x^Ts = \frac{1}{n} \sum_{i=1}^n x_i s_i,
\end{equation}
satisfies $\mu \le \epsilon$ for some tolerance $\epsilon>0$.

\section{Convergence and Complexity in Optimization.} \label{sec:complexity}

Useful optimization algorithms produce a sequence of iterates that tend toward solutions of the problem at hand.
Analysis of algorithms in optimization is concerned with proving results about the nature of these convergence properties.
For decades, most researchers working on nonlinear optimization focused on two kinds of convergence results.
\begin{itemize}
    \item Asymptotic (``global'') convergence: What is the ultimate behavior of the algorithm as more and more iterates are generated?  
    Does the sequence of iterates converge to a point satisfying certain optimality conditions (or at least have accumulation  points that satisfy optimality conditions)? Or stop at such a point?
    \item Local rate-of-convergence: Once the iterate sequence reaches a certain neighborhood of a point satisfying certain optimality conditions, how rapidly does it converge to that point?
\end{itemize}
In this paper, we focus mostly on {\em non-asymptotic} analysis, in which we ``globalize'' the local analysis and seek to say something about the rate of convergence of the algorithm from its initial point.
Although it has a long history, this kind of analysis has become much more prominent in recent years.
Non-asymptotic analysis of a particular algorithm identifies some scalar quantity  $\tau_k$ associated with each iteration or each oracle (see below), such that $\tau_k \to 0$ as we approach a point satisfying desired optimality conditions.  
For the unconstrained problem \cref{eq:f}, $\tau_k$ might be $f(x^k)-f^*$, $\| \nabla f(x^k) \|$, or $\mbox{\rm dist} (x^k,\cS)$, where $x^k$ is the $k$th iterate of the algorithm.
Depending on the algorithm, convergence of $\tau_k$ to zero can occur at arithmetic, ``sublinear'' rates, such as $\tau_k \le A/k$, $\tau_k \le A/\sqrt{k}$, or $\tau_k \le A/k^2$ for some constant $A>0$, which depends on the problem and often on the starting point.
It can occur at ``linear'' (sometimes called ``exponential'' or ``geometric") rates such as $\tau_k \le A (1-\phi)^k$, for some constants $A>0$ and $\phi \in (0,1)$.
It could occur at ``superlinear'' rates, for which $\tau_{k+1}/\tau_k \to 0$ as $k \to \infty$.

\paragraph{Iteration Complexity.}
Non-asymptotic properties provide bounds on the amount of computation required to satisfy the optimality conditions to some specified accuracy $\epsilon$.
The bounds on the convergence indicator $\tau_k$ as a function of $k$ can be converted into bounds on the number of iterations $k$ needed to ensure that $\tau_k \le \epsilon$ for any tolerance $\epsilon$ via a simple inversion.
For example, if $\tau_k \le A/k$ for some $A>0$, then we have $\tau_k \le \epsilon$ whenever $k \ge A/\epsilon = O(\epsilon^{-1})$. For linear convergence, when $\tau_k \le A (1-\phi)^k$, we have $\tau_k \le \epsilon$ when $k \ge (1/\phi) |\log (\epsilon/A)| = O(|\log \epsilon|)$.

Complexity analysis of this type is sometimes referred to as {\em iteration complexity}, since ``$k$" refers to the iteration index of the algorithm.
But this approach may not be immediately relatable to the conventional understanding of complexity in computer science, which concerns bounds on the amount of computation or storage required to solve problems in a class. 
The issues are that (a) within some algorithms, different iterations may require significantly different amounts of computation  --- and sometimes the amount of computation is hard to bound; (b) different instances within a problem class may require vastly different amounts of time to execute a single iteration. 
(It is easier to relate classical notions of complexity in linear programming than in general smooth minimization.)

\paragraph{Operation Complexity.}
Formally, we assume the Blum-Shub-Smale (BSS) model of complexity \cite{smale2000algorithms,blum2012complexity} in which the primitive objects are real numbers, and each arithmetic operations $+$, $-$, $\times$, $\div$ (as well as comparisons $\le$, $\ge$, and $=$) are each assumed to be a single unit of computation. 
(This model is closer to practical computation with floating-point numbers than the ``bit-complexity'' model, which assumes that problem data is rational and takes the unit of computation to be a bitwise operation.)
Essentially, the BSS model requires us to bound the number of real-number {\em arithmetic operations} required.
In interior-point LP algorithms, each iteration requires the assembly and solution of a system of linear equations of dimension $O(n)$, so we can bound the number of operations required per iteration by $O(n^3)$. 
if the number of iterations required to reduce the quantity \cref{eq:mu} below $\epsilon$ is $O(\sqrt{n}| \log \epsilon |)$, say, then the arithmetic complexity bound is $O(n^{3.5} | \log \epsilon|)$. 
(In fact, slightly tighter bounds can be obtained by modifying the algorithm to take advantage of the similarity between the linear equations that are solved at successive iterations; see \Cref{sec:lp}.)

\paragraph{Oracle Complexity.}
For nonlinear problems such as \cref{eq:f}, the {\em oracle complexity} model of Nemirovski and Yudin \cite{NemY83} is widely used  to bound the amount of computation required by a certain algorithm on a given class of problems.
For a given algorithm, this method defines some ``unit of information" (oracle) that the algorithm requires about the function $f$, with many queries made during the course of the algorithm. 
For one algorithm (e.g. a gradient descent method), a suitable oracle might by the value of $(f(x),\nabla f(x))$ at a given point $x$.
For another method (e.g. a stochastic gradient method), the appropriate oracle might be an unbiased random estimate of $\nabla f(x)$ with a certain known variance, for a given $x$.
The ``oracle complexity'' of a method is a bound on the number of oracles required to obtain an $\epsilon$-accurate solution.
When each iteration of an algorithm requires one oracle (or a fixed number of oracles), there is an obvious relationship between oracle complexity and iteration complexity.
This relationship is not always so tight; some algorithms (for example, a gradient-descent method with backtracking) might require a variable number of oracles per iteration.

\paragraph{Lower Bounds and Optimal Algorithms.}
Most complexity analyses are concerned with upper bounds on the relevant measure of computation. 
But there has also been much interest in {\em lower bounds}, which are usually defined in terms of both a class of algorithms and a class of problems.
The lower bound is usually derived from a worst-case instance in the problem class, making a  claim of the following form:
\begin{quote}
    There is an instance in the problem class for which any algorithm in the algorithm class will require at least $N(\epsilon)$ oracle calls to find an $\epsilon$-accurate solution.
\end{quote}
Algorithms for which the lower bound is within a constant multiple (not depending on $\epsilon$) of the upper bound is called an {\em optimal algorithm}. 
Possibly the most famous instance of an optimal method is Nesterov's accelerated gradient method, which is optimal in the algorithm class in which $x^k-x^0$ lies in the span of all gradients encountered by the algorithm between iteration $0$ (initial point $x^0$) and iteration $k$, and the problem class is the class of all convex functions whose gradients are Lipschitz continuous.

\paragraph{Gaps between Theory and Practice.}
For some algorithms, the (upper-bound) complexity of an algorithm reflects its typical behavior on a numerical instance. 
In other cases, the practical behavior is usually much better than the worst-case complexity analysis would suggest. 
This can happen for several reasons.
\begin{itemize}
    \item The algorithm design and analysis is based on assumptions about the problem that are much too conservative for most instances, or are tight on only a small fraction of the search space navigated by the algorithm. 
    For instance, an assumption that the gradient of $f$ in \cref{eq:f} is Lipschitz continuous with some constant $L$ is common in gradient-based methods for this problem. But the local geometry of the problem might see the algorithm making much faster progress in some regions of the parameter space, where the problem has much nicer properties than the worst-case assumptions would suggest.
    \item The  instances in a problem class that are difficult for an algorithm might be extremely rare. 
    Such is the case in linear programming, where instances for which the simplex method performs poorly are almost never encountered in practical or randomly generated instances.
    \item Algorithms may contain adaptive mechanisms (for example, line searches or trust-region strategies) that allow them to exploit variations in the properties of problems across the parameter space. 
    These gains may be reflected in practice but the variability of properties of the problem may be hard to express in an abstract way that lends itself to theoretical analysis.
    \item The problems of interest in a given class actually belong to a subclass with properties that set them apart from the wider class and make them easier to solve. 
    One example is the ``benign nonconvexity'' phenomenon, which has been encountered in many problems (especially from machine learning) over the past 10 years, where global minima of certain nonconvex objectives are usually found easily \cite{Sun21}, despite global minimization of general nonconvex objectives being intractable.
\end{itemize}
For these and other reasons, practitioners are generally not advised to use complexity bounds as the sole basis for choosing which algorithm to apply to a given problem.
Such bounds may be a consideration, but computational experience on similar problems is a more reliable guide.
Nevertheless, the process of understanding the gaps and perhaps narrowing or closing them --- by refining the complexity analysis, refining problem classes to identify subclasses that admit tighter analysis, or discovering algorithmic features or even new algorithms with better theoretical performance 
--- deepens our understanding of both the algorithms and problem classes.
Sometimes the insights so gained can percolate into wide practical use, as in the interior-point revolution in linear programming.

\paragraph{Historical Perspectives.}

The  influential early reference for complexity of algorithms for smooth convex nonlinear optimization is the book of Nemirovski and Yudin \cite{NemY83}. 
This book, which appeared in Russian in 1979, described the key concepts of problem class and algorithm class, oracle complexity, and lower bounds on complexity of methods, among other things.
Nesterov's text \cite{Nes04}, in its revised and expanded version\cite{nesterov2018lectures}, is a key modern reference.
It describes upper and lower complexity bounds for gradient methods and Newton methods for smooth functions and for subgradient and cutting-plane methods for nonsmooth functions, along with a description of barrier methods based on self-concordant functions and their complexity.
The recent book of Cartis, Gould, and Toint~\cite{cartis2022evaluation} describes complexity for smooth nonlinear optimization, both constrained and unconstrained, with a focus on nonconvex problems with higher-order optimality guarantees.


\section{Linear Programming.} \label{sec:lp}

Linear programming \cref{eq:lp} and the simplex method were developed to address wartime logistical planning in the 1940s \cite{dant:histo}.
LP quickly became a subject of fascination to the fledgling operations research and theoretical computer science communities, and has remained so to the present day.
The restriction to linear models and  ``certain" data ($A$, $b$, and $c$) were recognized as being unrealistic in practice, but the approximation is often good enough for the model to be useful. 
Moreover, LP is a foundational model for linear {\em integer} programming, in which some components of the variable $x$ in \cref{eq:lp} are restricted to have integer or binary ($0/1$) values. 
The most successful approaches for solving integer LPs (branch-and-bound, branch-and-cut) require the solution of many (closely related) LPs. 

In this section we outline the progress that has been made in algorithms and complexity theory for linear programming, and discuss briefly the practical significance of these various developments.
A fine summary of the history up to 2000 appears in \cite{Todd2002}; our discussion below makes note of progress in the years since then.

\subsection{Simplex Method.}

The simplex method starts with the recognition that the feasible set $\cF := \{ x \in \R^n \, : \, Ax=b, \; x \ge 0\}$ is a polyhedron in $\R^n$ --- a convex set with flat boundary surfaces and extreme points called {\em vertices}. 
Each vertex is the unique minimizer of some linear function $v^Tx$ over $\cF$, for some $v \in \R^n$. 
Moreover, when solutions to the LP \cref{eq:lp} exist, at least one such solution is a vertex, so nothing is lost by searching for a solution only among the vertices of $\cF$.
Each vertex can be characterized algebraically by forming a partition of the indices $\{1,2,\dotsc,n\}$ into sets $B$ and $N$ (corresponding to ``basic'' and ``nonbasic'' variables, respectively) for which $|B| = m$ and $|N| = n-m$. 
The basic set $B$ defines those components of $x$ which may be nonzero at the vertex, whereas all components $x_i$, $i \in N$ are fixed at zero.
In other words, the vertex satisfies the algebraic conditions 
\begin{equation} \label{eq:vertex}
    Ax=b, \quad \mbox{$x_i \ge 0$ for all $i \in B$}, \quad \mbox{$x_i = 0$ for all $i \in N$}.
\end{equation}
This characterization may not be unique; there may exist multiple partitions $B \cup N$ that define the same vertex, a phenomenon known as degeneracy.
Moreover, not every subset $B  \subset \{1,2,\dotsc,n\}$ of size $m$ necessarily yields a vertex, as there may be no point $x$ that satisfies the conditions \eqref{eq:vertex} for this $B$.
Nevertheless, it is easy to construct examples of feasible sets $\cF$ for which the number of vertices is exponential in dimension $n$.

Each iteration of simplex method proceeds by moving from one vertex to an adjacent vertex, obtained by swapping a single element $i \in \{1,2,\dotsc,n\}$ between sets $B$ and $N$, a process called ``pivoting."
The choice of element $i$ can usually be done in a way that both maintains feasibility of the resulting vertex and causes a reduction in the objective $c^T x$.
(Sometimes, the step is {\em degenerate}, in that $B$ and $N$ change but the vertex $x$ stays the same.)
When no adjacent vertex can be found with a lower value of the objective, the current vertex must be a solution, and the algorithm terminates.

At most iterations, there are multiple choices of the index $i$ that yield a decrease in $c^Tx$, and different variants of the simplex method have different {\em pivot rules} to identify some of these possibilities and choose one pivot from among them. 
From the invention of the method, the number of iterations / pivots required for the simplex method to converge on practical instances was observed to be at most a modest multiple of $n$.
Could a worst-case bound --- say a bound that is {\em polynomial} in $n$ --- be proved to hold for any problem of the form \eqref{eq:lp}, for the simplex method with a certain pivot rule? 
The question is still open, but important progress has been made.
In 1972, Klee and Minty  \cite{KleM72} provided a famous example for which Dantzig's original pivot rule requires a number of pivots exponential in $n$.
Numerous other pivot rules have been investigated, including some that are used in software implementations and some that are randomized (that is, they choose the pivot in some randomized way from among the set of possibilities that improve the objective).
In each case, instances of \eqref{eq:lp} have been constructed where the number of pivots is exponential or subexponential in $n$.\footnote{A subexponential lower bound is one in which the number of pivots is bounded below by $c_1 \exp (c_2 n^c)$, for constants $c_1>0$, $c_2>0$, and $c \in (0,1)$.}
(For the randomized methods, these complexities are for the {\em expected} number of pivots.)
See \cite[Section~1]{disser2023exponential} for a discussion of the current state of the art. 

\paragraph{Average-case Analyses.}
There have been a number of efforts in the 70s and 80s to find a bound on the average number of iterates required by some variant of the simplex method over some distribution of instances of \eqref{eq:lp}.
Borgwardt was the first to obtain a polynomial bound, and refined his approach over a number of years; see for example \cite{Bor87,Borgwardt1999}. 
Generally, the matrix $A$ in \cref{eq:lp} is chosen to have its columns dense in $\R^m$, distributed independent, identically, and symmetrically under rotations. 
The vectors $b$ and $c$ are chosen to make the problem feasible, in a well-defined way.
Taking the expectation over this problem distribution, the number of simplex pivots is $O(n^{1/(m-1)} m^2)$.
This type of analysis, while an important contribution, has little relevance to LP instances arising in practice, which are not distributed in the same way as the random instances in these studies. 
(For one thing, the matrix $A$ tends to be sparse in practical LPs.)

\paragraph{Smoothed Analysis.}
A breakthrough in theoretical understanding of the simplex method came in 2004 with the {\em smoothed analysis} of Spielman and Teng~\cite{SpeT04}.
Their result works with the dual form \cref{eq:lp.dual2} and a particular variant of the simplex method, known as the shadow-vertex simplex method.
Given a particular instance of \cref{eq:lp.dual2} defined by the triplet $(A,c,b)$, they produce perturbed instances obtained by adding Gaussian random variables to all components of $A$ and $b$, with mean zero and standard deviation $\sigma \max_{i=1,2,\dotsc,n} \, \| [A_{\cdot,i},b_i] \|$, where $[A_{\cdot,i},b_i]$ is the vector in $\R^{m+1}$ consisting of the $i$th column of $A$ and the $i$th element of $b$.
Taking the {\em expected} number of steps of their simplex variant over this set of perturbed random instances, they show that it is bounded by a quantity polynomial in $m$, $n$, and $1/\min(\sigma,\sigma_0)$, for some $\sigma_0>0$ \cite[Theorem~5.1]{SpeT04}. 
This result shows, roughly speaking, that even instances of \cref{eq:lp.dual2} for which simplex has poor behavior can likely be modified to an instance with polynomial-time behavior by adding small random perturbations to all the elements of $A$ and $b$.
Significant progress has been made since the original paper in improving the order of the polynomial dependence of the number of pivots; see for example \cite{Vershynin2009,bach2025optimalsmoothedanalysissimplex}. 
In particular, the dependence on standard deviation of the noise $\sigma$ has been reduced from $\sigma^{-30}$ in \cite{SpeT04} to $\sigma^{-1/2}$ in \cite{bach2025optimalsmoothedanalysissimplex}, with the latter also proving a matching high-probability lower bounding instance for which the number of steps is also a multiple of  $\sigma^{-1/2}$.

\subsection{Polynomial-Time Methods.}
The lack of progress in proving polynomial worst-case bounds for the simplex method spurred a search over several decades for alternative methods with this property. These efforts began to bear fruit in the late 1970s, and soon led to new approaches that have had a major impact on both theory and computations in LP.

\paragraph{Ellipsoid Method.}
The ellipsoid method was developed by Yudin and Nemirovskii \cite{YudinNemirovski1976} for convex nonlinear programming.
In 1979, Khachiyan~\cite{Kha79} achieved a breakthrough when he showed that an adaptation of this approach to LP converged in polynomial time --- the first polynomial-time algorithm for LP.
The method starts with an ellipsoid large enough to contain all the solutions, and each subsequent iteration generates successively smaller ellipsoids, all of which still contain the  set of solutions.
The method is usually described with respect to the dual problem \eqref{eq:lp.dual2}, which has variable $\lambda \in \R^m$ and a feasible set defined by $n$ linear inequalities.
At iteration $k$, it tests whether the center of the current ellipsoid in $\R^m$ is a solution to this problem.
If not, it is easy to identify a hyperplane that cuts away more than half the volume of the ellipsoid but has all the solution points in the remaining volume. 
A new ellipsoid is then constructed which contains this remaining volume, with a new center and a smaller volume. 
The key to the approach is that the volume of the containing ellipsoid  decreases by a fraction of approximately $\exp(-\frac{1}{2m})$ at each iteration. 
In ${O}(m^2 \log \epsilon)$ iterations, a solution can be identified to within an accuracy of approximately $\epsilon$.
The amount of arithmetic required at each iteration $k$ is $O(mn)$, since we need to form the product $A^T \lambda^k$, where $\lambda^k$ is the center of the ellipsoid at iteration $k$.
(An additional $O(m^2)$ operations are needed to update the center and matrix defining the ellipsoid.)

Although the ellipsoid method was an intriguing approach of great theoretical significance, its impact on practical computations with LP was negligible.
Unlike for the simplex method, the worst-case complexity turned out to be typical, and much slower than the practical performance of simplex. 
Moreover, there were issues of stability, requiring more accurate computations to be performed as the iterations progressed.

\paragraph{Karmarkar's Projective Algorithm.}
Karmarkar's announcement in 1984 \cite{Kar84} of an algorithm for LP with polynomial-time convergence and good practical performance, caused a sensation.
His algorithm considered a special form of the constraints in the LP \eqref{eq:lp}, namely,
\begin{equation}
    \label{eq:lp.k}
    \min_{x \in \R^n} \, c^Tx \;\; \mbox{subject to} \;\; Ax=0, \;\; \bfone^Tx = 1, \;\; x \ge 0,
\end{equation}
where $\bfone$ is the vector in $\R^n$ whose components are all $1$. 
Moreover, the point $(1/n) \bfone$ was assumed to be feasible for \eqref{eq:lp.k}.\footnote{A reformulation can be applied to the standard form \eqref{eq:lp} to obtain this special form of the constraints.}
The iterates $x^k$ of Karmarkar's method are all {\em strictly} feasible, that is, $Ax^k=0$, $\bfone^T x^k=1$, $x^k > 0$ ---  
a property that gave rise to the term ``interior-point method''.
At iteration $k$,  we define a new variable $z \in \R^n$ by $z \leftarrow (X^k)^{-1} x$, where $X^k = \diag (x^k_1,x^k_2, \dotsc, x^k_n)$, and consider a slightly modified LP:
\begin{equation} \label{eq:k2}
\min_{z \in \R^n} \, (X^k c)^Tz \;\; \mbox{subject to} \;\; (A X^k)z=0, \;\; \bfone^T z =1, \;\; z \ge 0.
\end{equation}
Note that $z = (1/n) \bfone$ is strictly feasible for this problem, as for \eqref{eq:lp.k}.
We take a step by projecting the cost vector $(X^k c)$ into the subspace defined by $(AX^k)v=0$, $\bfone^Tv=0$, and taking a step from $(1/n) \bfone$ along the negative of this projected direction, choosing the distance to maintain strict positivity of all components.
The new iterate $x^{k+1}$ is obtained by applying the diagonal scaling matrix $X^k$ to the resulting vector.

This {\em projective} method (so named because of its use of the projection of the rescaled cost vector) was shown in \cite{Kar84} to require $O(n \log \epsilon)$ iterations to reduce the objective by a factor of $\epsilon$, with each iteration requiring about $O(n^{2.5})$ arithmetic operations.
In contrast to the ellipsoid method, however, the practical performance of Karmarkar's method was much better than the worst-case analysis would suggest, especially when some flexibility was allowed in the choice of steplength.

The idea of diagonal rescaling using the current iterate was quickly discovered to have intriguing connections to earlier work.
In 1967, Dikin~\cite{Dik67} had proposed a similar algorithm using this rescaling.
After Karmarkar's announcement, \cite{GilMSTW86} showed that the search direction in the projective algorithm could be related to one obtained from a log-barrier approach, in which the constraints $x \ge 0$ are accounted for by subtracting a  term $\mu \sum_{i=1}^n \ln x_i$ (for some $\mu>0$) from the objective.

\paragraph{Interior-Point Methods: Theory and Practice in the  ``Classical" Period.}
Although the initial claims about Karmarkar's method being competitive computationally with simplex were not borne out, his announcement led to an explosion of activity in interior-point methods that reverberated for over a decade. 
Many polynomial-time approaches were proposed, during this ``classical" period,  for solving  the primal \eqref{eq:lp} or dual \eqref{eq:lp.dual} formulations, or solving both simultaneously by finding a triple $(x,\lambda,s)$ satisfying the optimality conditions \eqref{eq:lp.opt}.
Methods of the latter type, known  as {\em primal-dual methods}, proved to be particularly fruitful as an area for development. 
The two major classes of methods in this are primal-dual potential reduction methods (proposed by Tanabe, Todd, and Ye \cite{Tan87,TodY90}) and path-following methods. 
The latter came in many flavors, proposed and studied by many authors; a brief outline of the history can be found in \cite[Chapter~5]{Wri97}. 
Iteration complexity bounds for almost all these methods ranged in a narrow band from $O(n^{1/2} \log \epsilon)$ and $O(n^2 \log \epsilon)$ iterations.

An important work of this period, which mostly adopts the primal rather than primal-dual viewpoint, was the book of Nesterov and Nemirovskii~\cite{NesN93}.
This work influenced theory and practice for convex problems including LP but extending further.
It  replaces the constraint $x \in \cF$ where $\cF$ is a closed convex set with a term $\mu \phi(x;\cF)$ added to the objective, where $\mu>0$ is a barrier parameter (gradually reduced to zero during the computation) and $\phi$ is a {\em barrier function} whose domain is the relative interior of $\cF$ with the property that $\phi(x;\cF) \to \infty$ as $x$ approaches the boundary of $\cF$. 
The barrier function satisfies an additional property of {\em self-concordance}, which (roughly speaking) allows its third derivatives to be bounded in terms of a $3/2$ power of the second derivatives,  as in the function $-\log t$ for $t \in \R$.
Armed with this property, they construct a generic algorithm based on a decreasing sequence of values of the barrier parameter $\{\mu_k\}$, with Newton steps taken on the sum of the objective and weighted barrier function for each $\mu_k$. 
The complexity of this approach is similar to the primal-dual LP methods discussed above, with a polynomial dependence on dimension $n$.
Barrier functions with the self-concordant property can be constructed explicitly for several convex optimization problems, including LP, convex quadratic programming, second-order cone programming, and semidefinite programming.
Computationally, these approaches are usually slower than primal-dual algorithms, but they provided rigor and inspiration for the approaches implemented in such codes as   SDPT3~\cite{sdpt3}, SeDuMi~\cite{Stu99}, and the commercial solver MOSEK~\cite{mosek}.

Parallel to the theoretical developments in path-following algorithms for LP, practical methods based of this type were being tested on test problems batteries, which were curated to be representative of real-world problems.
(The "netlib" test set (\url{https://www.netlib.org/lp/}) was particularly influential.)
Good practical interior-point codes started to appear by the end of the 1980s, with both commercial and public-domain codes appearing throughout the 1990s.
After the first few years, and up to the present, the algorithm underlying almost all interior-point software for LP has been Mehrotra's predictor-corrector primal-dual approach \cite{Meh92a}, which is a path-following method with clever heuristics to select certain critical parameters.
Another key development was that the quality of commercial implementations of the simplex method improved vastly following the appearance of interior-point competition.
The combination of new algorithmic ideas (in the simplex, interior-point, and sparse linear algebra spaces), better implementations, and more powerful computers has led to an enormous increase in the size of LPs that can be solved practically by current software and by many orders of magnitude reduction in runtime.

Early interior-point software tended to follow the form of one or other theoretical method but with liberties taken in the choice of parameters, including steplengths. 
Quickly, however, the algorithms that were being implemented became quite similar to the ones that were being analyzed theoretically. 
Gaps remain between the worst-case complexity for primal-dual methods (see above) and the practical experience, where the number of iterations grows only very slowly with $n$. 
(These gaps are trivial by comparison with the corresponding gaps for the simplex method.)
An interesting point is that the most successful practical primal-dual approach is closest to the methods with $O(n^2 \log \epsilon)$ complexity; methods with the slightly better $O(n^{1/2} \log \epsilon)$ bound are somewhat slower in practice.

\paragraph{Primal-Dual Interior Point Methods.}
We give a little more detail on path-following primal-dual interior-point methods, which in addition to having a strong complexity theory have given rise to the most successful interior-point implementations.
These methods formulate the LP optimality conditions \eqref{eq:lp.opt} as the following system of bound constrained, (mildly) nonlinear equations:
\begin{equation}
    F(x,\lambda,s) := \left[\begin{matrix}
        A^T \lambda + s -c \\ Ax-b \\ XS \bfone 
    \end{matrix}\right] = \left[\begin{matrix}
        0 \\ 0 \\ 0
    \end{matrix}\right], \quad (x,s) \ge 0,
\end{equation}
where $X = \diag(x_1,x_2, \dotsc,x_n)$ and $S=\diag(s_1,s_2,\dotsc,s_n)$, so that $XS \bfone$ is the vector with components $x_is_i$, $i=1,2,\dotsc,n$.
Note that the system $F(x,\lambda,s)=$ is ``square;" the number of equations and number of unknown is the same.
Path-following steps start by defining a central path, which is the set of strictly feasible points for which the products $x_i s_i$, $i=1,2,\dotsc,n$ are all identical. That is, the central path is the set of points $\{(x(\mu), \lambda(\mu), s(\mu))$ parametrized by a scalar $\mu>0$ that satisfies
\begin{equation} \label{eq:cp}
    F(x,\lambda,s;\mu) = \left[\begin{matrix}
        A^T \lambda + s -c \\ Ax-b \\ XS \bfone - \mu \bfone
    \end{matrix}\right] = \left[\begin{matrix}
        0 \\ 0 \\ 0
    \end{matrix}\right], \quad (x,s) > 0,
\end{equation}
(Under mild conditions, including the full row rank of the matrix $A$ that we assumed at the start, this system has a unique solution for each $\mu>0$.)
Path-following methods enforce a condition slightly stronger than strict feasibility on each of their iterates, requiring essentially that the pairwise products $x_i s_i$, $i=1,2,\dotsc,n$ are not too different from each other. 
A set of this type is $\cN_{-\infty}(\gamma)$ for some parameter $\gamma \in (0,1)$, defined as follows:
\[
\cN_{-\infty}(\gamma) := \{ (x,\lambda,s) \, : \, Ax=b, \; A^T \lambda+s=c, \; (x,s)>0, \; x_i s_i \ge \gamma (x^Ts)/n, \; i=1,2,\dotsc,n \}.
\]
Given any $(x,\lambda,s) \in \cN_{-\infty}(\gamma)$, we can apply Newton's method to the equality constraints in \eqref{eq:cp}, for some choice of $\mu$, to obtain the following linear system for the step $(\Dx,\Dl,\Ds)$:
\begin{equation} \label{eq:lpf1}
    \nabla F(x,\lambda,s;\mu) 
    \left[ \begin{matrix}
    \Dx \\ \Dl \\ \Ds 
    \end{matrix} \right] = -F(x,\lambda,s;\mu) \quad
    \Leftrightarrow \quad
    \left[ \begin{matrix}
    0 & A^T & I \\ A & 0 & 0 \\ S & 0 & X 
    \end{matrix} \right]
    \left[ \begin{matrix}
    \Dx \\ \Dl \\ \Ds 
    \end{matrix} \right] =
    - \left[\begin{matrix}
        0 \\ 0 \\ XS \bfone - \mu \bfone
    \end{matrix}\right].
\end{equation}
The following algorithm computes a series of steps of this type, decreasing the ``target" value of $\mu$ at each iteration, and choosing the steplength to be as large as possible, subject to staying in the neighborhood $\cN_{-\infty}(\gamma)$ (and not exceeding $1$).

\begin{algorithm} \label{alg:lpf}  \caption{{\bf LPF:} Long-Step Path-Following Algorithm for \eqref{eq:lp}}
\begin{algorithmic}
    \STATE Given $\gamma$, $\sigmin$, and $\sigmax$ with $\gamma \in (0,1)$ and $0 < \sigmin < \sigmax <1$, and $(x^0,\lambda^0,s^0) \in \cN_{-\infty}(\gamma)$; 
    \FOR{$k=0,1,2,\dotsc$}
    \STATE Define $\mu_k := (x^k)^T s^k/n$;
    \STATE Choose $\sigma_k \in [\sigmin,\sigmax]$;
    \STATE Solve  \eqref{eq:lpf1} with $\mu = \sigma_k \mu_k$ and $(x,\lambda,s) = (x^k,\lambda^k,s^k)$ to obtain $(\Dx^k,\Dl^k,\Ds^k)$;
    \STATE Choose $\alpha_k$ to be the largest value of $\alpha$ in $[0,1]$ such that 
    \[
    (x^k,\lambda^k,s^k) + \alpha (\Dx^k,\Dl^k,\Ds^k) \in \cN_{-\infty} (\gamma);
    \]
    \STATE Set $(x^{k+1},\lambda^{k+1},s^{k+1}) = (x^k,\lambda^k,s^k) + \alpha_k (\Dx^k,\Dl^k,\Ds^k)$;
    \ENDFOR
\end{algorithmic}
\end{algorithm}

Remarkably, this simple algorithm can be shown to find a point $(x^T,\lambda^T,s^T)$ satisfying \eqref{eq:lp.opt} in an $\epsilon$-approximate sense (that is, the first three conditions in \eqref{eq:lp.opt} are satisfied and in place of the last condition we have $(x^T)^Ts^T/n\le \epsilon$) in $T=O(n |\log \epsilon|)$ iterations, provided that the initial point $(x^0,\lambda^0,s^0)$ satisfies $(x^0)^Ts^0/n \le 1/\epsilon^{\kappa}$ for some $\kappa>0$ independent of $n$ and $\epsilon$.
This powerful result requires only elementary mathematics and can be proved from scratch in a single lecture of a graduate-level optimization course (see \cite[pp.~96--100]{Wri97}).
Moreover, for reasonable choices of the parameters $\sigmin$, $\sigmax$, and $\gamma$, it works quite well in practice.
Even so, there is a tension between the theory and the practice.
The constant factor in the $K=O(n |\log \epsilon|)$ iteration bound depends on $\gamma$ through a term $(1+\gamma)/(\gamma (1-\gamma))$, which is smallest when $\gamma \approx 1/2$. 
Better practical results are obtained for smaller values of $\gamma$, such as $\gamma=10^{-3}$, which allows the pairwise products $x_i s_i$ to vary widely.

Mehrotra's ``predictor-corrector'' method \cite{Meh92a}, still the basis of most practical interior-point codes, broadly tracks the LPF algorithm above but makes several modifications that enhance its practical performance (but not the complexity bounds).
It allows iterates to be infeasible with respect to the conditions $Ax=b$, $A^T\lambda+s=c$ (but reduces the infeasibility to zero as the algorithm proceeds); makes an adaptive choice of the parameter $\sigma_k$ in Algorithm~\ref{alg:lpf}, which typically results in more aggressive progress ($\sigma_k$ closer to zero) as the algorithm nears the solution; and incorporates a second-order correction into the step calculation by modifying the right-hand side of \eqref{eq:lpf1}. 
For details, see \cite[Chapter~10]{Wri97}.

Primal-dual algorithms, including the path-following approaches described here, can be extended straightforwardly, beyond LP to wider classes of problems such as convex quadratic problems and monotone linear complementarity problems. 
Extensions to second-order cone programming and semidefinite programming are also possible, though these are somewhat more complicated; see for example~\cite{Tod01a,TutTT03}.
In all these contexts, similar to LP, there are path-following methods with polynomial complexity, and the practical method of choice is a (nontrivial) extension of Mehrotra's predictor-corrector method.

\begin{figure} \label{fig:yintat}
    \centering\includegraphics[width=5in]{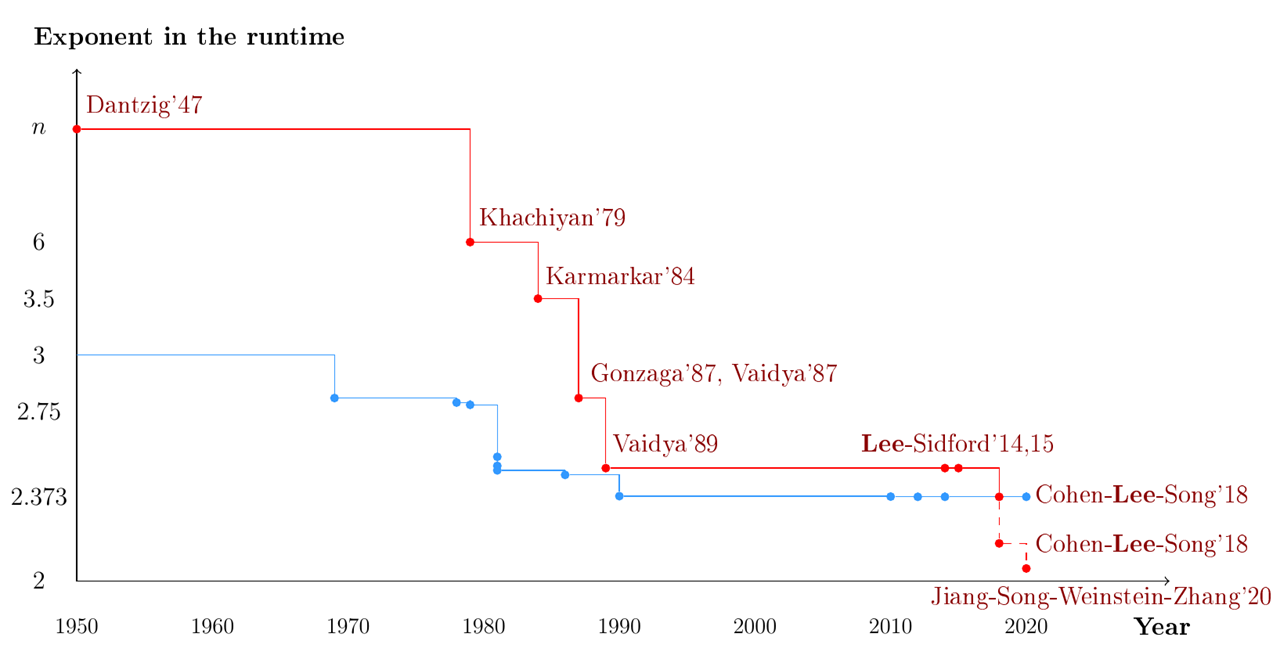}
    \caption{Progression of complexity bounds in matrix multiplication (blue) and linear programming (red). The final dashed piece of the red line shows the complexities of the methods of  \cite{cohen2021solving,van2020deterministic} and \cite{jiang2021faster}, respectively, that would apply if the exponent of matrix multiplication were shown to be $2$. (Credit: Yin-Tat Lee.)}
\end{figure}

 \paragraph{Recent Developments in Complexity Analysis of Interior-Point Methods.} 
For many years after the ``classical'' period of interior-point linear programming research ended in the mid-1990s, there was little research on improving the complexity bounds further.
A burst of new activity started in the late 2010s, resulting in refinements of these bounds, for both general linear programs \cref{eq:lp} and various special cases such as network flows and graph matching.
We will focus on the algorithms for general LP, described  first as a stochastic algorithm in  \cite{cohen2021solving} and later modified by \cite{van2020deterministic} (which presents a slightly simpler deterministic algorithm)  and \cite{jiang2021faster}.
All are based on a primal-dual path-following method that solves Newton-like linear systems of the form \cref{eq:lpf1} at each iteration, restated here in general form as follows:
\begin{equation}
    \label{eq:lpf3}
    \left[ \begin{matrix}
    0 & A^T & I \\ A & 0 & 0 \\ S & 0 & X 
    \end{matrix} \right]
    \left[ \begin{matrix}
    \Dx \\ \Dl \\ \Ds 
    \end{matrix} \right] =
    \left[\begin{matrix}
        0 \\ 0 \\ r
    \end{matrix}\right],
\end{equation}
where $r=-(XS \bfone - \mu \bfone)$ in \cref{eq:lpf1}. 
The components $\Dx$ and $\Ds$  of the solution of \cref{eq:lpf3} can be written as:
\begin{equation}
    \label{eq:lpf4}
    \Dx = X (XS)^{-1/2} (I-P) (XS)^{-1/2}  r, \quad \Ds = S(XS)^{-1/2} P (XS)^{-1/2} r, 
\end{equation}
where the projection matrix $P$ is defined as 
\begin{equation}
    \label{eq:lpf5}
    P = W^{1/2} A^T (A W A^T)^{-1} A W^{1/2}, \quad W = XS^{-1}.
\end{equation}
(Note that since $S$ and $X$ are positive diagonal matrices, $W$ is also positive diagonal, as are any products and powers involving $X$, $S$, and $W$.)
The technique used in \cite{cohen2021solving,van2020deterministic,jiang2021faster} is to take path-following steps, reducing the quantity $\mu_k=(x^k)^Ts^k/n$ by a factor of roughly $(1-c/\sqrt{n})$ at each iteration (for some $c>0$ independent of $n$), so that an $\epsilon$-accurate solution can be found in $O(\sqrt{n} \log \epsilon)$ steps.
But the number of operations required to solve an {\em approximate} version of \cref{eq:lpf3} on each iteration is reduced far below the naive estimate of $O(n^3)$. 
The key observation is that the diagonals of the matrix $w$ (namely, $x_i/s_i$, $i=1,2,\dotsc,n$) change only slowly in a relative sense from iteration to iteration, so they maintain an approximate version of $W$ --- denoted $\bar{W}$ ---  in which the diagonal elements are updated only if they have changed significantly since their last update.
Since the number of such updates is small on each iteration, there is only a low-rank change to the matrix $A\bar{W}A^T$, so the inverse of this matrix, and ultimately the approximate solution of \cref{eq:lpf3} at the next iteration, can be calculated cheaply. 
Given the current best-known exponent $\omega = 2.371339$  \cite{alman2024asymmetryyieldsfastermatrix} of matrix multiplication\footnote{``Exponent $\omega$'' means that the complexity of multiplying two $n \times n$ matrices is $O(n^{\omega})$.}, the overall computational cost of these methods is $\tilde{O}(n^{\omega} \log \epsilon)$ (where the notation $\tilde{O}$ indicates that poly-logarithmic factors in $n$ are omitted).
This is a remarkable fact: In principle, the cost of linear programming is within a poly-log factor of the cost of matrix multiplication!
Were the exponent of matrix multiplication  ultimately shown to be $2$ (as is conjectured, although progress in reducing it has slowed considerably), the complexity of the algorithms in \cite{cohen2021solving,van2020deterministic} would become $\tilde{O}(n^{2+1/6})$ while for the algorithm of \cite{jiang2021faster} it would become $O(n^{2+1/18})$. 
\Cref{fig:yintat} shows the progression of best-known complexities for matrix multiplication (blue line) and linear programming (red line), with the rightmost two points on the red line showing the complexities of the methods just discussed under the assumption of $\omega=2$.

\subsection{First-Order Methods.}
Since 2021, there has been renewed interest in formulating LP as a saddle-point problem and applying first-order methods, with many enhancements, to solve it \cite{applegate2021practical,applegate2025pdlp,lu2024restarted}.
Such ideas had been tried in the 1980s without much success on practical problems, but several factors gave rise to  renewed interest.
First, the enormous size of some modern LP applications make the storage and computational requirements of  matrix factorizations, required by both simplex and interior-point methods, difficult even for modern computing hardware.
Second, numerous enhancements to the basic first-order algorithms have been made over the past several years, and have been essential to them becoming competitive.
Third, since their main computational requirement is to form matrix-vector products involving the constraint matrix $A$ and its transpose, they are much more amenable to computation using GPUs (graphical processing units) than either simplex or interior-point methods. 

The primal-dual pair \cref{eq:lp}, \cref{eq:lp.dual2} can be written equivalently as the following min-max problem:
\begin{equation}
    \label{eq:lp.minmax}
    \min_{x \ge 0} \, \max_{\lambda} \, c^Tx - \lambda^T Ax + b^T \lambda.
\end{equation}
Given a starting point $(x^0,\lambda^0)$, the PDHG algorithm \cite{chambolle2011first} applied to this problem alternates between updating $x$ and updating $\lambda$ according to the following scheme, for $k=0,1,2,\dotsc$:
\begin{equation} \label{eq:lp.pdhg}
    x^{k+1} = \left[ x^k - \tau (c-A^T\lambda^k) \right]_+, \quad 
    \lambda^{k+1} = \lambda^k + \sigma (b-A(2x^{k+1}-x^k)),
\end{equation}
for certain positive choices of stepsizes $\tau$ and $\sigma$, where $[ \cdot]_+$ denotes projection onto the nonnegative orthant.
A Halpern variant described in \cite{lu2024restarted} defines the new iterate $(x^{k+1},\lambda^{k+1})$ as a weighted average of the iterate from  \eqref{eq:lp.pdhg} and the initial point $(x^0,\lambda^0)$. Recent experience indicates some theoretical and computational advantages for the Halpern approach \cite{lu2025cupdlpx}.

Classically, under certain conditions on the choice of stepsizes, and subject to existence of a solution satisfying \eqref{eq:lp.opt}, the running average of the iterates generated by \eqref{eq:lp.pdhg} converges to an $\epsilon$-accurate solution in $O(1/\epsilon)$ iterations. 
Since LP has a property called {\em sharpness} (essentially that a certain measure of primal-dual optimality grows linearly with distance to the solution set), this rate can be improved to linear, that is, $O(|\log \epsilon|)$ iterations for an $\epsilon$-accurate solution. 
However, the constant multiplier of $|\log \epsilon|$ is problem-dependent and almost always extremely large --- it is related to the Hoffman constant $H(A,b,c)$ of the system of linear equalities and inequalities that defines the solution, namely,
\begin{equation} \label{eq:lp.opt2}
Ax=b, \quad A^T \lambda \le c, \quad  x \ge 0, \quad c^Tx - b^T \lambda = 0
\end{equation}
(equivalent to \cref{eq:lp.opt}).
The theory can be refined to take account of the fact that as the algorithm approaches a solution, the active set (essentially, those components of $x$ that are zero at the solution) is identified, so the behavior  of the method is governed by a reduced version of \cref{eq:lp.opt2} that is also homogeneous, so it typically has a much  smaller Hoffman constant. 
The practical performance tends to match this theory; acceleration is noted near the solution.
Many enhancements have been applied to these methods to make them competitive in practice with state-of-the-art solvers, for 
 very large LPs. 
 One important enhancement is {\em restarting}, which accelerates the linear rate of convergence both in theory (reducing the iteration bound from $O(H(A,b,c)^2 |\log\epsilon|)$ to $O(H(A,b,c) |\log\epsilon|)$) and in practice.


Recalling the average-case analysis of the simplex method discussed above, we note a recent result \cite{xiong2025high} about LPs drawn from a certain distribution (that includes the elements of $A$ being i.i.d. sub-Gaussian with zero mean and unit variance).
With probability $1-\delta$, the number of iterations required by restarted PDHG is $\tilde{O}(P(m,n) |\log \epsilon|/\delta)$, where $P(m,n)$ is polynomial in dimensions $m$ and $n$, and $\tilde{O}$ hides log factors. 
(Here $\delta$ is small and positive, but not exponentially small in terms of $m$ and $n-m$.)

\section{Unconstrained Optimization.} \label{sec:f}

We consider now the unconstrained smooth minimization problem \cref{eq:f}, which requires minimization of a function $f:\R^n \to \R$.
We will assume that $f$ is as smooth as it needs to be for the algorithm at hand.
Typical assumptions are that $f$ is differentiable with Lipschitz continuous gradients, that is, there is a constant $L$ such that 
\begin{equation} \label{eq:glip}
  \| \nabla f(x) - \nabla f(y) \| \le L \|x-y\|, \quad \mbox{for all $x,y \in \R^n$},
\end{equation}
or (stronger) that $f$ is twice differentiable with Lipschitz continuous Hessian, that is, there is a constant $M$ such that
\begin{equation} \label{eq:hlip}
  \| \nabla^2 f(x) - \nabla^2 f(y) \| \le M \|x-y\|, \quad \mbox{for all $x,y \in \R^n$},
\end{equation}

Naturally, the theoretical and practical difficulty of minimizing $f$ varies greatly with the properties of $f$, even among problems of the same dimension $n$.
Convexity provides a fundamental way to taxonomize these problems.
A smooth function is {\em convex} if for all $x,s \in \R^n$ we have
\begin{equation}
  \label{eq:convex}
  f(x+s) \ge f(x) + \nabla f(x)^Ts,
\end{equation}
and {\em strongly convex} if for some $\mu>0$ (termed the ``modulus of convexity''), we have
\begin{equation}
  \label{eq:sconvex}
f(x+s) \ge f(x) + \nabla f(x)^Ts + \frac{\mu}{2} s^Ts.
\end{equation}
If the function is twice continuously differentiable, these conditions are equivalent to $\nabla^2 f(x) \succeq 0$ and $\nabla^2 f(x) \succeq \mu I$, respectively.

As we noted earlier, the first-order conditions for $x^*$ to be a
local minimizer of $f$ are that $\nabla f(x^*)=0$; second-order
necessary conditions add the condition $\nabla^2 f(x^*) \succeq 0$;
while second-order sufficient conditions add the condition $\nabla^2
f(x^*) \succ 0$. In fact, the latter conditions ensure that $x^*$
has a strictly smaller function value than any other points in a
neighborhood, and that there are no other local minimizers in this
neighborhood.  Note that for {\em convex} $f$, first-order points
automatically satisfy second-order necessary conditions.
(They satisfy  second-order
{\em sufficient} conditions for strongly convex $f$.)
In fact, points
satisfying these conditions are not just local solutions but {\em
  global} solutions; their optimality holds not just in a neighborhood
but in the entire space $\R^n$.

When $f$ is not convex, there may in general be infinitely many local
solutions, with many different function values. (Consider for example
the scalar function $f(x) = x^2/10 + \sin (\pi x)$. Usually we are
most interested in a global solution --- one for which $f(x^*)\le
f(x)$ for all $x$ --- but the problem of finding this is intractable
for general smooth nonconvex $f$.  Algorithms for minimizing $f$ guarantee
convergence to a point satisfying first-order or second-order
conditions, possibly with a complexity guarantee for finding points that satisfy these conditions {\em approximately}.
An interesting development of the past decade has been the discovery of many classes of functions, mostly from machine learning and data science applications, that are nonconvex but which (for various reasons) global solutions can be found tractably, with complexity guarantees. We discuss this point further below.

The simplest kind of nonlinear smooth function is a quadratic, which has the form
\begin{equation} \label{eq:quad}
  f(x) = \frac12 x^THx + g^Tx,
\end{equation}
where $H$ is an $n \times n$ symmetric matrix.
For this function we have $\nabla f(x) = Hx+g$ and $\nabla^2 f(x) =H$.
This function is unbounded below when $H$ has one or more negative eigenvalues; we can make $f$ arbitrarily negative by moving far enough in the direction of an eigenvector corresponding to one of the negative eigenvalues. Thus the problem of minimizing $f$ has no solution.
When $H \succeq 0$ but there exist {\em zero} eigenvalues, the problem has a solution only if $g^Ts=0$ for all $s \in \R^n$ with $s^THs=0$. (In this case, there is an affine space of points that are all global minimizers of $f$.)
When $H$ is positive definite, the problem is strongly convex, with a unique minimizer $x^* = -H^{-1} g$.

Convex quadratics serve as a kind of sanity check on gradient-based
methods for unconstrained minimization. Such methods correspond to
iterative methods for solving the square symmetric linear system
$Hx=-g$ (corresponding to $\nabla f(x)=0$) and can be analyzed
with the tools of linear algebra. However,  algorithmic strategies such as multistep descent and conjugate gradients developed for convex quadratics often need to be modified significantly for case of more general nonlinear functions, and their convergence properties may carry over only in a weaker form, or not at all.

The {\em finite-sum} form of $f$ has become enormously important in machine learning:
\begin{equation} \label{eq:finite-sum}
  f(x)= \frac{1}{N} \sum_{i=1}^N f_i(x),
\end{equation}
where each $f_i:\R^n \to \R$ is a smooth function.
In the machine learning context, $x$ contains the parameters of a model (for example, the ``weights'' in  a neural network) while each $f_i$  measures the accuracy of the model's prediction on the $i$th element of training data.
Often, $N$ is extremely large, making even the evaluation of $f$ impossible to perform in reasonable computing time --- let alone its gradient and Hessian.
The techniques used to minimize finite-sum $f$ are strongly motivated by the gradient algorithms of this section, but differ in that they use unbiased estimates of the gradient $\nabla f$ at each iteration, based on an average over a relatively small sample of the component-wise  gradients $\nabla f_i(x)$, $i=1,2,\dotsc,N$.

\subsection{Gradient Methods.} \label{sec:grad}

We consider first methods that take steps of the form
\begin{equation} \label{eq:gd}
  x^{k+1} = x^k - \alpha_k \nabla f(x^k), \quad k=0,1,2,\dotsc,
\end{equation}
where $\alpha_k>0$ is the {\em steplength} and $x^0$ is the chosen
starting point.  These methods are often called {\em steepest-descent}
or {\em gradient} methods.  (The term ``steepest'' refers to the direction giving the greatest reduction in $f$ per
distance moved, for small values of this distance, according to the
first-order Taylor series approximation.)

\paragraph{Descent Methods.}
Traditionally, $\alpha_k$ in \eqref{eq:gd} is chosen so as to reduce
$f$ --- not just infinitesimally (which would not justify the cost of
calculating $\nabla f(x^k)$) but by a significant amount.  An ``exact''
line search chooses $\alpha_k = \arg\min_{\alpha >0} \, f(x^k+ \alpha
s^k)$, where $s^k:= -\nabla f(x^k)$, but finding such an $\alpha$ is
in itself a (scalar) minimization problem, usually requiring multiple
evaluations of $f$ and the directional derivative $\nabla
f(\cdot)^Ts^k$. In practice, we settle for a crude approximation to
this optimal steplenth that satisfies certain conditions to guarantee
convergence yet is efficient to compute. One such condition is a
``sufficient decrease'' condition, which has the form
\begin{equation} \label{eq:suff}
f(x^k+\alpha s^k) \le f(x^k) + c_1 \alpha \nabla f(x^k)^T s^k, \quad \mbox{for some $c_1 \in (0,\tfrac12)$}.
\end{equation}
In fact, all sufficiently small values of $\alpha$ satisfy this
condition, as the first-order Taylor series expansion yields
$f(x^k+\alpha s^k) \approx f(x^k) + \alpha \nabla f(x^k)^T s^k$, and $
\nabla f(x^k)^T s^k = -\| \nabla f(x^k) \|^2$ for the negative
gradient direction $s^k = -\nabla f(x^k)$. (The condition
\cref{eq:suff} is usually coupled to another condition to ensure that
$\alpha$ is not too short; see, for example \cite[Chapter~3]{NocW06}).

We can get a more precise bound on the reduction in $f$ by assuming
that $f$ has Lipschitz continuous gradients \cref{eq:glip}, in which
  case we have a global upper bound on the function:
\begin{equation} \label{eq:qbound}
    f(x^k+\alpha s^k) \le f(x^k) + \alpha \nabla f(x^k)^Ts^k + \alpha^2 \frac{L}{2} \| s^k \|^2.
\end{equation}
For $s^k =-\nabla f(x^k)$, the value of $\alpha$ that minimizes the
right-hand side is $\alpha = 1/L$, and the condition \cref{eq:suff}
will hold for all $\alpha \in (0,2(1-c_1)/L)$. 

To obtain complexity bounds under various assumptions on $f$, we use
the constant-steplength choice $\alpha_k \equiv 1/L$ along with $s^k =
-\nabla f(x^k)$, so the algorithm \cref{eq:gd} becomes
\begin{equation} \label{eq:gdc}
  x^{k+1} = x^k - \frac{1}{L} \nabla f(x^k), \quad k=0,1,2,\dotsc,
\end{equation}
(Similar complexity bounds can be derived for methods that make more
practical choices of $\alpha_k$, whose computational performance is
usually --- but not always --- better.) By substituting this choice of
$\alpha_k$ into \cref{eq:qbound}, and using \cref{eq:gd}, we obtain
\begin{equation} \label{eq:grad1}
  f(x^{k+1}) \le f(x^k) - \frac{1}{2L} \| \nabla f(x^k) \|^2.
\end{equation}
Suppose that $f$ is lower-bounded, that is, there is $\bar{f}$ such
that $f(x) \ge \bar{f}$ for all $x$. By rearranging \eqref{eq:grad1}
and taking the average over $k=0,1,\dotsc,T-1$ for some $T$, we obtain
\begin{equation} \label{eq:grad2}
  \frac{1}{2L} \frac{1}{T} \sum_{k=0}^{T-1} \| \nabla f(x^k) \|^2 \le \frac{1}{T} \sum_{k=0}^{T-1} \left( f(x^k) - f(x^{k+1}) \right) \le \frac{f(x^0) - \bar{f}}{T}.
  \end{equation}
By rearranging and using the fact that the minimum  is smaller than the average, we have
\begin{equation} \label{eq:grad3}
  \min_{k=0,1,\dotsc,T-1} \, \| \nabla f(x^k) \| \le \frac{\sqrt{2L(f(x^0)-\bar{f})}}{\sqrt{T}}.
  \end{equation}
In summary, if we run the algorithm \cref{eq:gdc} on a smooth
function $f$ for at least $T = O(\epsilon^2)$ iterations, then at
least one of these iterates $k=0,1,\dotsc,T$ we will encounter an
iterate $x^k$ for which the approximate first-order condition $\|
\nabla f(x^k) \| \le \epsilon$ holds.

For convex $f$, a similar analysis gives a faster rate of convergence,
for the same algorithm \cref{eq:gdc}. Supposing that there exists a
minimizer $x^*$ of $f$, we have from convexity \cref{eq:convex} that
$f(x^*) \ge f(x^k) + \nabla f(x^k)^T(x^*-x^k)$. By substituting into
\cref{eq:qbound} and using \cref{eq:gdc} and rearranging, we have
\[
f(x^{k+1}) \le f(x^*) + \nabla f(x^k)^T(x^k-x^*) - \frac{1}{2L} \| \nabla f(x^k) \|^2 =
f(x^*) + \frac{L}{2} \left( \| x^k-x^* \|^2 - \| x^{k+1} - x^* \|^2 \right).
\]
By subtracting $f(x^*)$ from both sides, averaging over iterates
$k=0,1,\dotsc,T-1$, telescoping the sum and using $f(x^T)-f(x^*)  \ge 0$, we obtain
\[
\frac{1}{T} \sum_{k=0}^{T-1} (f(x^{k+1})-f(x^*)) \le \frac{L}{2T} \| x^0-x^* \|^2.
\]
Since $f$ decreases at each iteration, we can bound the average
suboptimality on the left-hand side below by $f(x^T)-f(x^*)$, to obtain
\begin{equation} \label{eq:comp.conv}
  f(x^T)-f(x^*) \le \frac{L}{2T} \| x^0-x^* \|^2.
\end{equation}
Thus, we need $T=O(\epsilon^{-1})$ iterations to find a point whose
function value is within $\epsilon$ of the optimal function value.
(\cite{pmlr-v97-lee19e} improves the rate \cref{eq:comp.conv} to $o(1/T)$.)

For the case in which $f$ is smooth and {\em strongly} convex
\cref{eq:sconvex} with modulus of convexity $\mu$, the minimizer $x^*$
is unique, and a slightly longer but still elementary analysis (see for example \cite[Chapter~3]{WrightRecht22}) yields that
\begin{equation} \label{eq:comp.sconv}
  f(x^T) - f(x^*)\le \left( 1-{1}/{\kappa} \right)^T (f(x^0)-f(x^*)),  \quad
\mbox{where $\kappa:= L/\mu$,}
\end{equation}
so that a point with $f(x^T)-f(x^*)$ can be guaranteed after at most
$T = \kappa \log ((f(x^0)-f(x^*))/\epsilon)$ iterations.\footnote{This bound can be approximately halved if the fixed steplength $\alpha_k \equiv 2/(L+\mu)$ is used in place of $1/L$.}
The quantity
$\kappa$ measures the conditioning of the problem. When $f$ is
quadratic, it measures the ratio of the largest to smallest
eigenvalues of $H$.

Note that for the convex and strongly convex cases, $x^*$ is a global
minimizer of the problem, whereas for the general case \cref{eq:grad3}
guarantees only that the {\em first-order conditions} $\nabla f(x^*)$
will be satisfied approximately.

The strongly convex case has only a logarithmic dependence on
$\epsilon$, as opposed to the arithmetic dependence on $\epsilon^{-1}$
exhibited for general convex functions and nonconvex functions, but
can still be slow when the conditioning is poor.
The paper \cite{necoara2019linear} catalogs a variety of conditions weaker than
strong convexity under which linear convergence similar to
\eqref{eq:comp.sconv} is still obtained.  The best known of these is
the Polyak-Lojasiewicz (PL) condition \cite{KarNS16a} which requires
that $\tfrac12 \| \nabla f(x) \|^2 \ge \mu(f(x) -f(x^*))$ for some
$\mu>0$ (for some solution $x^*$).

When the steplength $\alpha_k$ in \cref{eq:gd} is chosen by some
adaptive line-search scheme, the worst-case theoretical bounds
essentially match those just proved for the fixed-step case. This is
because the key inequality \cref{eq:grad1} holds, with the factor
$1/(2L)$ replaced by some other positive constant.  The practical
performance of such adaptive algorithms is usually better, but there are classes of functions for which the worst-case
behavior is typical. 

Given that even exactly-minimizing choices of $\alpha$ cannot yield
better worst-case results than those for the constant $1/L$ step, we
are tempted to conclude that there can be no other schemes using the
search direction $s^k =-\nabla f(x^k)$ that yield better complexities.
Remarkably, this is not true! The key is to not demand reduction of
$f$ at every step, but rather to take longer steps on some iterations
over which $f$ {\em increases}, an algorithmic property called {\em
  nonmonotonicity}. The idea is that such steps might move to a new
iterate that allows larger reduction to take place at the next step
(or some later step), so that over a span of multiple consecutive
steps, we may see a larger overall improvement in $f$ than over the
equivalent number of iterations of the ``greedy'' approach motivated by the single-iteration perspective.

\paragraph{Nonmonotone Methods.}
The nonmonotone algorithm proposed in 1988 by Barzilai and Borwein
\cite{BarB88}, drew inspiration from the Newton and quasi-Newton
methods described later in this section. In the basic Newton method,
the search direction is $ -\nabla^2 f(x^k)^{-1} \nabla
f(x^k)$. \cite{BarB88} propose to form a one-parameter approximation
$\alpha_k I$ to the inverse Hessian, by choosing $\alpha_k$ to imitate
the action of the Hessian over the step just taken, from $x^{k-1}$ to
$x^k$. Defining $s^k = x^k-x^{k-1}$ and $y^k = \nabla f(x^k) - \nabla
f(x^{k-1})$, Taylor's theorem suggests that $s^k \approx \nabla^2
f(x^k)^{-1} y^k$. Consequently, they choose $\alpha_k$ to satisfy $s^k
\approx \alpha_k y^k$ in the least-squares sense, leading to the
explicit formula $\alpha_k = (s^k)^T y^k / (y^k)^Ty^k$. This choice of
$\alpha_k$ in \eqref{eq:gd} gives an algorithm that can increase
measures of convergence alarmingly on some iterations, but sees
dramatic improvements on other iterations. \cite{BarB88} give a short, novel analysis of a 2-variable convex quadratic that
explains these observations nicely. However, there is no analysis that
proves better worst-case rates than those described above for the
constant-stepsize methods on general functions.

The Barzilai-Borwein approach has been extended in many ways in the
decades since, with new variants and modifications to improve global
convergence properties and to handle bound constraints and nonsmooth
regularization. The practical performance on nonlinear problems is
similar to nonlinear conjugate gradient (a monotone method, discussed
below) and requires a similar amount of storage. The originality of
the approach, its unusual theoretical properties, and its good
practical performance on a number of problem types (see
\cite{WriNF08a} for an example) have made it consistently popular.


Recently, another approach \cite{altschuler2025acceleration} has been proposed for algorithms of the
form \cref{eq:gd} that are nonmonotone and that have better worst-case
complexities than those of the basic gradient descent method, on
convex functions. 
This work considers sequences of steplengths that yield the optimal improvement in distance to the solution $\| x-x^* \|$ over multiple consecutive steps.
An antecedent of this approach that applies to strictly convex quadratics \cref{eq:quad} was proposed by Young in 1953 \cite{young1953richardson}, who showed that the optimal steplengths in \cref{eq:gd} over a sequence of $m$ consecutive steps correspond to the inverse roots of an order-$m$ Chebyshev polynomial. 
The improvement is strictly better than $m$ steps with any fixed steplength.
The general idea of multistep optimization carries over to general smooth, strictly convex functions --- but in a way that it is not at all straightforward.
The steplength schedule for quadratics is provably bad for the more general case \cite[Chapter~8]{altschuler2018greed}.
Instead, \cite{altschuler2025acceleration} describes a recursive construction of ``silver stepsize schedules" for values of $m$ that are powers of $2$.
When plotted, these schedules are fractal and approximately periodic, with period of about $\kappa^{\log_{1+\sqrt{2}} 2} \approx \kappa^{.7864}$, where $\kappa=L/\mu$ is the condition number for $f$.
Occasionally, steplengths much longer than the fixed step $1/L$ discussed above are taken.
Construction of the schedules is explicit but quite intricate; the convergence theory uses some of the analytical machinery developed during the past 15 years for accelerated gradient methods (see below).
Ultimately, it is shown \cite[Theorem~1.1]{altschuler2023acceleration} that a final error of $\|x_T - x^* \| \le \epsilon$ can be attained in about $T=O(\kappa^{.7864} \log (1/\epsilon))$ iterations. 
Note that the dependence on $\kappa$ is significantly better than for the fixed-step method (where the exponent $.7864$ is replaced by $1$) but not as fast as the accelerated gradient methods discussed below, which have complexity bounds $O(\kappa^{1/2} \log (1/\epsilon))$.
As we will see, the latter methods take steps in directions other than the negative gradient $-\nabla f(x^k)$, although the search direction in each is a linear combination of gradients evaluated at points visited by the algorithm up to and including the current iterate.

This approach is extended to smooth convex (but not necessarily strongly convex) functions in \cite{altschuler2023acceleration}, with complexity bound $O(\epsilon^{-.7864})$, improving over the $O(\epsilon^{-1})$ bound \cref{eq:comp.conv} for the fixed-steplength method.
Another recent multistep approach for functions in this class \cite{grimmer2024provably} leads to stepsize schedules that also feature occasional very long steps, but does not improve the dependence of worst-case complexity on $\epsilon$.

Practical implications of the methods in \cite{altschuler2025acceleration,altschuler2023acceleration,grimmer2024provably} are unclear. 
Among methods using only gradients and making use of known properties of $f$ (such as $L$ and $\mu$), they are dominated in theory by the accelerated gradient methods discussed below, and they lack the adaptivity features of \cite{BarB88}.
However, they provide important support for the notions that the multistep perspective can lead to better theoretical results, that nonmonotonicity can be tolerated, and that occasional long steps can lead to better overall behavior.

\paragraph{Convergence to Second-Order Points.}

When applied to nonconvex functions $f$, most algorithms of the form \cref{eq:gd} are guaranteed at best to converge to stationary points, for which $\nabla f(x^*)=0$. 
Such  points may not be local minima, for example, when  the Hessian $\nabla^2 f(x^*)$ has negative eigenvalues, so that $s$ is a direction of descent for $f$ and $x^*$ if $s^T \nabla^2 f(x^*) s<0$. (Such a point is called a {\em strict saddle} point.)
\cite{Lee16h} showed that gradient descent \cref{eq:gd} with constant steplength $\alpha_k \equiv \alpha \in (0,1/L]$ almost certainly does not converge to strict saddle points when $x^0$ is chosen randomly.
(The result is proved by using a stable manifold theorem from the dynamical systems literature; see for example \cite{Shub87}.)
Extension of this result to the heavy-ball method (described below) is shown in \cite{OneW18a}, along with analysis of the rates of escape of gradient and accelerated gradient methods from the neighborhood of strict saddles.

A modification to the gradient method \cref{eq:gd} was proposed in \cite{jin2017escape}, in which small perturbations are added occasionally to $x^k$ when $\| \nabla f(x^k) \|$ falls below a small threshold.
After such a perturbation is added, their method waits to see if subsequent gradient descent steps escape the neighborhood of $x^k$, by observing how $f(x)$ and $\| \nabla f(x) \|$ change on subsequent iterations.
If $f$ does not decrease and $\| \nabla f(x) \|$ remains small, there is evidence that $x^k$ is not a strict saddle point, and in fact approximately satisfies the second-order conditions $\nabla f(x)=0$, $\nabla^2 f(x) \succeq 0$, in the sense that
\begin{equation} \label{eq:approx2o}
  \| \nabla f(x^k) \| \le \epsg, \quad \nabla^2 f(x^k) \succeq -\epsH I,
\end{equation}
for small tolerances $\epsg$ and $\epsH$.
In fact, \cite{jin2017escape}  show that the number of gradient evaluations required is bounded by $O(\epsilon^{-2})$ when  $\epsg = O(\epsilon)$ and $\epsH = O(\sqrt{\epsilon})$ in \cref{eq:approx2o}.
A subsequent work by these authors \cite{JinNJ17a} improves the complexity to $O(\epsilon^{-7/4})$, approximately matching the oracle complexity of Newton-based methods described below, by using an accelerated gradient method and a mechanism to exploit negative curvature in the Hessian.

\subsection{Accelerated Gradient Methods.} \label{sec:ag}

We turn now to methods that require evaluation of gradients at each iteration (sometimes one gradient and sometimes several), but which achieve faster non-asymptotic rates than the methods of the previous subsection, at least in some settings. 
A feature  common to all these methods is that the direction of the step at each iteration is generally not $-\nabla f(x^k)$, as in \cref{eq:gd}, but rather a linear  combination of all the gradients encountered so far during the algorithm.

\paragraph{Conjugate Gradient Methods.}
The notion that faster progress toward minimizing $f$ in \cref{eq:f} could be made not just by stepping along the gradient at the current point but by aggregating information from gradients calculated earlier in the algorithm, dates to at least the 1950s, with the advent of the conjugate gradient (CG) method of Hestenes and Steifel \cite{Hestenes}.
This method applies to minimization of strongly convex quadratics (equivalently, solution of linear systems with symmetric positive definite coefficient matrices).
Starting from a point $x^0$ with gradient $r^0 := \nabla f(x^0) = Hx^0+g$, it generates iterates $x^k$ with the remarkable property that each $x^{k+1}$ minimizes $f$ in the affine space $x^0 + \cK_k(H;r^0)$, where 
\[
\cK_k(H;r^0) := \mbox{span}\{r^0, Hr^0, H^2 r^0, \dotsc, H^kr^0 \} = \mbox{span} \{ \nabla f(x^0), \nabla f(x^1), \dotsc, \nabla f(x^k) \}
\]
is known as the {\em $k$th Krylov subspace} for $H$ generated by $r^0$.
In general, the dimension of this space grows by $1$ at each iteration, until a minimizer is found --- in at most $n$ iterations.
For search directions $p^k$,  CG sets $x^{k+1} = x^k + \alpha_k p^k$ for $\alpha_k$ that exactly solves  $\alpha_k = \arg\min_{\alpha >0} \, f(x^k+ \alpha p^k)$,
The set of directions $\{ p^0, p^1, \dotsc \}$ is constructed so that $p^0=r^0$, $\cK_k(H;r^0) = \mbox{span} \{ p^0,p^1, \dotsc, p^k \}$, and  the property of  {\em conjugacy} holds with respect to $H$, that is, $(p^i)^T H p^j =0$ when $i \neq j$.
Each  $p^k$ is generated at iteration $k$, according to the formula 
\begin{equation} \label{eq:cg.beta}
p^{k} = -\nabla f(x^k) + \beta_k p^{k-1},
\end{equation}
where $\beta_k>0$ is calculated by a simple formula that ensures that the properties described above hold.
Descriptions of the method and its properties can be found in \cite[Chapter~10]{GolV83} and \cite[Section~5.1]{NocW06}.

Besides finite termination in at most $n$ steps, this {\em linear CG} method has powerful theoretical properties.
Convergence is measured in terms of the weighted  norm $\| x-x^* \|_H^2 = (x-x^*)^TH(x-x^*)$, with an overall linear rate of 
\begin{equation} \label{eq:cg}
\| x^k-x^* \|_H^2 \le  \left( 1- \frac{2}{\sqrt{\kappa}+1} \right)^k \|x^0-x^*\|_H^2,
\end{equation}
where $\kappa = L/\mu$ is the condition number of $H$ (the ratio of largest to smallest eigenvalue).
It follows that $\| x^T - x^* \|_H \le \epsilon$ can be achieved in $O(\sqrt{\kappa} \, |\log (\epsilon/\|x^0-x^*\|_H)|)$ iterations.
More refined analyses of convergence based on the distribution of eigenvalues of $H$ are available, showing that faster convergence occurs when the eigenvalues are clustered.

Linear CG is a workhorse algorithm of scientific computing, being widely used to solve PDEs. 
For practical implementations, good preconditioning (a transformation of $Hx=-g$ into an equivalent system $(C^{-T} H C^{-1}) \hat{x} = -C^{-1} g$ for some nonsingular $C$) is essential to make the spectrum of the Hessian more conducive to fast convergence. 
The conjugacy property of the vectors $\{ p^0,p^1, \dots \}$ degrades in practice because of roundoff error. 
However, the method has maintained its primacy among iterative techniques for solving $Hx=-g$ for several decades.

Unfortunately, adaptations of the linear CG algorithm to nonlinear $f$ do not retain the powerful properties of the quadratic case, either in theory or practice --- even for strongly convex $f$.
The earliest nonlinear CG algorithm, proposed in 1964 by Fletcher and Reeves \cite{FleR64} is the most straightforward extension, choosing $\alpha_k$ as an approximate minimizer of the line search along $p^k$, and replacing the formula for the parameter $\beta_k$ in \eqref{eq:cg.beta}  by a sensible generalization. 
Like many nonlinear CG methods, the Fletcher-Reeves (FR) method reduces to the original CG method when $f$ is strongly convex quadratic and the line search for $\alpha_k$ is exact.
No convergence rates are known, but even for nonconvex $f$, the gradients $\nabla f(x^k)$ generated by this method accumulate at zero.
A later method proposed in 1969 by Polak and Ribi\`ere  \cite{PolR69} and Polyak \cite{Polyak1969} uses a different formula for the parameter $\beta_k$.
An enhanced version of this method (known as ``PR$+$'' or ``PRP$+$'') restarts by setting $\beta_k=0$ (equivalently, $p^k = -\nabla f(x^k)$) to ensure that $p^k$ is always a descent direction for $f$.
This approach generally performs better than the FR approach in practice, with similar convergence properties \cite{gilbert1992global}.
\cite{Polyak1969} proved a complexity bound of $O(\kappa^3 \log \epsilon)$ for strongly convex $f$ with condition number $\kappa$ for the unmodified PR algorithm.

The variant of nonlinear CG due to Hager and Zhang \cite{HagZ05}, has theoretical and practical properties that are possibly the best for this class of methods.
Its formula  for $\beta_k$  ensures that every $p^k$ is a descent direction while not explicitly restarting, and in fact allowing $\beta_k<0$ in some circumstances.
Theoretical properties include convergence to a minimizer for strongly convex $f$ (but at rates unknown) and accumulation of gradients at zero for general $f$, as for the earlier methods.
The same authors have written an  excellent review of nonlinear conjugate gradient methods \cite{HagZ06b}.

The recent paper \cite{das2024nonlinear} contains new  results on both worst-case complexity for the PRP and FR  variants on strongly convex nonlinear $f$, assuming that the line search for $\alpha_k$ is exact in each case.
The results use sophisticated computer-aided proofs, in which quadratically constrained quadratic programs (QCQP) are solved to find worst-case $f$ (among smooth functions with modulus of complexity $\mu>0$ and gradient Lipschitz constant $L$) for which the improvement in optimality gap $f(x^k)-f(x^*)$ over a single iteration (or multiple iterations) is as weak as possible.
Essentially, neither PRP nor FR has better complexity than gradient descent with exact line search (which converges to $f(x^T)-f(x^*)$ in $T = O(\kappa \log \epsilon)$ iterations), except in the case of a slight improvement for PRP when $\kappa$ is large. 
The bounds do not approach the $O(\sqrt{\kappa} \log \epsilon)$ bounds observed for linear CG and for the accelerated gradient methods discussed below.

Karimi and Vavasis \cite{KarV21a} hybridize the Hager-Zhang nonlinear CG method with Nesterov's accelerated gradient method (described below), in a way that achieves the best-known rates on convex and strongly convex nonlinear functions, while recovering the finite termination properties of linear CG if $f$ is a strongly convex quadratic. 
Essentially, the approach uses CG but monitors a measure of progress derived from Nesterov's analysis of his accelerated gradient methods, and switches to that method when insufficiently rapid progress is made. 
When the function appears to be near-quadratic, it switches back to CG.
This method combines the best available complexity results with good numerical performance.


In \Cref{sec:newton}, we will see another use of CG in the optimization context, as a solver for the linear system of equations to be solved at each iteration of Newton's method. 
This way of using CG allows us to exploit the strong properties of both linear CG and Newton's method, to obtain algorithms that are practical and for which theoretical guarantees are available.

\paragraph{Heavy-Ball.}
The heavy-ball method introduced by Polyak in 1964 \cite{Pol64a} shares some features with conjugate gradient, but is simpler. Its steps have the form
\begin{equation} \label{eq:hb}
x^{k+1} = x^k - \alpha_k \nabla f(x^k) + \beta_k (x^{k}-x^{k-1}), \quad k=0,1,2,\dotsc,
\end{equation}
where $x^0$ is the starting point, $x^{-1}=x^0$,  and $\alpha_k>0$ and $\beta_k \ge 0$ are scalar parameters.
The relationship with CG is evident when we compare \cref{eq:hb} with \eqref{eq:cg.beta}.
In both cases, the step at each iteration is a linear combination of the negative gradient direction and the previous step, with the previous step itself being a linear combination of gradients encountered at previous iterates.
But whereas in the case of CG the relative weights of these two components and the overall scaling applied to the step has an intricate motivation, based entirely on the case of convex quadratic $f$, the motivation of heavy-ball is more elementary.
It can be viewed as steepest decent \cref{eq:gd}  with the addition of a ``momentum'' component $\beta (x^{k}-x^{k-1})$: We continue to move in the direction of the previous step, tweaked to incorporate the latest gradient information.
Momentum introduces a smoothing effect, recognizing that a weighted averaging of information across all gradients may provide a more reliable indication of the overall landscape of the function than just the latest gradient.
When the parameters $\alpha$ and $\beta$ are tuned appropriately, heavy-ball can indeed generate a smoother trajectory than the zig-zagging often associated with steepest descent (see, for example \cite[Figure~6]{Pol87}).

For the case of strongly convex quadratic $f$ \cref{eq:quad}, worst-case convergence rates for heavy-ball are similar to for linear CG.
Precisely, \cite[Section~3.2.1]{Pol87} shows that  for any $\delta>0$, and making the fixed-parameter choices $\alpha_k \equiv 4\mu/(\sqrt{\kappa}+1)^2$ and $\beta_k \equiv (\sqrt{\kappa}-1)^2/(\sqrt{\kappa}+1)^2$ in \cref{eq:hb}, there are constants $C_\delta$ and $K_\delta$ such that 
\begin{equation} \label{eq:hb.conv}
\| x^k-x^* \| \le C_\delta \| x^0-x^* \| (\rho^*+\delta)^k, \quad \mbox{for all $k>K_{\delta}$,}
\end{equation}
where $\rho^* = (\sqrt{\kappa}-1)/(\sqrt{\kappa}+1)$. Note that $\rho^* \in (1-2/\sqrt{\kappa}, 1-1/\sqrt{\kappa})$, so that, roughly speaking, the number of iterations $T$ required to ensure that $\|x^T -x^* \| \le \epsilon$ is $T = O(\sqrt{\kappa} | \log \epsilon|)$.

Although this result captures the global behavior of the heavy-ball algorithm for convex quadratic, and possibly also the ``local'' behavior of the algorithm once the iterates enter a neighborhood of the minimizer $x^*$, it does not extend to the case of general nonlinear strongly convex $f$. 
Indeed, Polyak's choice of parameters $\alpha_k$, $\beta_k$ does not even guarantee convergence for strongly convex nonlinear $f$ (\cite{LessardRechtPackard2016}; see  below).
Different choices of parameters can yield global linear convergence, but the best rate established so far leads to a complexity of $T = O(\kappa |\log\epsilon|)$ \cite{ShiDuJordanSu2022} --- no better than (unaccelerated) gradient descent --- and this bound appears to be tight \cite{goujaud2023provable}.

A recent variation on heavy ball \cite{weiaccelerated25} modifies \cref{eq:hb} as follows:
\begin{equation} \label{eq:hb.mod}
x^{k+1} = x^k - \alpha_k (2 \nabla f(x^k) - \nabla f(x^{k-1})) + \beta_k (x^{k}-x^{k-1}), \quad k=0,1,2,\dotsc,
\end{equation}
with $\alpha_k \equiv \mu/(\sqrt{\kappa}+1)^2$ and $\beta_k \equiv  1/(\sqrt{\kappa}+1)^2$. 
This method attains $f(x^T)-f(x^*) \le \epsilon$ in $T=O(\sqrt{\kappa} |\log\epsilon|)$ iterations, the same rate as is seen for CG and \cref{eq:hb} in the quadratic case. 
For general convex functions \cref{eq:convex}), a further modification of \cref{eq:hb.mod} requires $T = O(\epsilon^{-1/2})$ iterations for $f(x^T)-f(x^*) \le \epsilon$ --- an acceleration over the rate  seen for gradient descent on convex functions \cref{eq:comp.conv}.

\paragraph{Nesterov Acceleration.}
In 1983, Nesterov \cite{Nes83} proposed a method for general smooth convex functions with the accelerated rate of $T = O(\epsilon^{-1/2})$ for $f(x^T)-f(x^*) \le \epsilon$. 
A variant of the approach for strongly convex functions also has a rate that accelerates over steepest descent: $f(x^T)-f(x^*) \le \epsilon$ in $T=O(\sqrt{\kappa} |\log\epsilon|)$ iterations.
These rates are achieved with a seemingly minor variation of the heavy-ball iteration formula \cref{eq:hb}:
\begin{equation} \label{eq:nest}
x^{k+1} = x^k - \alpha_k \nabla f(x^k + \beta_k (x^{k}-x^{k-1})) + \beta_k (x^{k}-x^{k-1}), \quad k=0,1,2,\dotsc,
\end{equation}
the difference being that the gradient $\nabla f$ is evaluated at the extrapolated point $x^k + \beta_k (x^{k}-x^{k-1})$ rather than at $x^k$. 
By introducing an intermediate variable $y^k$, \cref{eq:nest} can be rewritten as:
\begin{equation}
    \label{eq:nest1}
    y^k = x^k + \beta_k (x^k-x^{k-1}), \quad x^{k+1} = y^k - \alpha_k \nabla f(y^k), \quad k=0,1,2,\dotsc.
\end{equation}
For the strongly convex case, the constant stepsizes $\alpha_k \equiv 1/L$ and $\beta_k \equiv (\sqrt{\kappa}-1)/(\sqrt{\kappa}+1)$ yield the accelerated rate. 
A variety of non-intuitive and highly technical proofs were proposed to demonstrate accelerated rates.
One approach is to define a Lyapunov function $V_k(x^k,x^{k-1})$ and show that it decreases linearly to zero at the desired rate.
A function of this type (see for example \cite[Section~4.3]{WrightRecht22}) is
\[
V(x^k,x^{k-1}) := f(x^k) - f(x^*) + \frac{L}{2} \| (x^k-x^*) - \rho^2 (x^{k-1}-x^*) \|_2^2,
\]
where $\rho^2 = (1-1/\sqrt{\kappa})$ for the given choices of $\alpha_k$ and $\beta_k$.
For the general convex case, the constant $\rho$ in this definition is replaced by a variable $\rho_k \in (0,1)$, which approaches $1$, while for the parameters in \cref{eq:nest} we can choose $\alpha_k \equiv 1/L$ while $\beta_k = \rho_k \rho_{k-1}^2$. which also approaches $1$. The Lyapunov analysis can be deployed to prove that $f(x^k)-f(x^*) = O(1/k^2)$, leading to the accelerated rate.

We discuss below alternative approaches to proving convergence and to designing new algorithms, that have been developed in the past 15 years.

\paragraph{Lower Bounds.}
The main focus on this paper is on {\em upper bounds} on complexity, with rates being identified for different algorithms on different classes of problems. 
It is natural to ask: Are these rates improvable?
An important tool for answering this question is {\em lower bounds} on the convergence rate.
For a class of algorithm that makes use of a certain kind of oracle, and a class of problems that these algorithms solved, is there some ``difficult" problem on which no algorithm in the class can achieve convergence faster than the lower-bounding rate?
If some algorithm in the class has an upper bound that is the same as the lower bound (to within a constant factor), it is called an {\em optimal algorithm}.

The classic lower bound, and the one most relevant to our present discussion, is due to Nesterov, and described in \cite[Theorem~2.7.1]{nesterov2018lectures} (see also \cite[Section~4.6]{WrightRecht22}).
The class of problems is smooth functions with Lipschitz smoothness $L$, and the class of algorithms are those for which $x^k-x^0$ (the aggregation of all steps since initialization at $x^0$) lies in the space spanned by all gradients evaluated up to iteration $k$.
(The algorithms discussed above all have this property.)
Nesterov's example is a convex quadratic \cref{eq:quad} with a tridiagonal Hessian, starting from $x^0=0$ and for which $x^*$ has all $n$ components nonzero. 
Each iterate $x^i$ has nonzeros only in its first $i$ components, so a lower bound on the distance $\| x^i-x^* \|$ is the norm of the last $n-i$ components of $x^*$.
For iterates $k=1,2,\dotsc, (n/2)-1$, we have $f(x^k)-f(x^*) \ge \frac{3L\|x^0-x^*\|^2}{32 (k+1)^2}$. 
Nesterov's algorithm is optimal for this class, as its upper-bounding rate is $O(1/k^2)$.
Note that the lower bound is valid only for  the first $n/2$ iterations. 
Limitations based on the dimension are standard in lower-bounding examples (see for example \cite{Drori2017}).

\paragraph{Accelerated Methods: Recent Developments for Convex $f$.}
The  theory of accelerated gradient methods has continued to make remarkable advances in recent years; see the monograph \cite{dAspremontScieurTaylor2021} for an excellent survey. 
The major development is that the processes of analyzing and designing accelerated methods like Nesterov's can be ``automated'' by expressing them as optimization (or saddle-point) problems which can be solved numerically and sometimes  even analytically. 
A breakthrough came around 2012 with the ``Performance Estimation Problem (PEP)'' technique Drori and Teboulle \cite{DroriTeboulle2014}. 
To obtain upper bounds for the rate of a given algorithm, they formulate an optimization problem to  {\em maximize} $f(x^T)-f(x^*)$ after $T$ steps of a given algorithm, the variable in this problem being the function $f$ and the constraints capturing the steps of the algorithm.
This formulation is not solvable as it stands, because $f$ lives in a function space (for example, ``$L$-smooth convex functions'') and is inherently infinite-dimension.
The key insight is to parametrize $f$ by variables that represent the iterates of the algorithm $x^i$, $i=0,1,\dotsc,T$,  the function values $f^i=f(x^i)$, $i=0,1,\dotsc,T$, and the gradients $g^i = \nabla f(x^i)$, $i=0,1,\dotsc,T$, and then introduce additional constraints to ensure that a function $f$ in the desired function space can be obtained by interpolating this function and gradient information.
The resulting problem becomes finite dimensional, and can be stated as a matrix optimization problem (the matrix variable being constructed from the gradients $g^i \in \R^n$, $i=0,1,\dotsc,T$).
After some ingenious relaxations, formulation tricks, and application of duality, a formulation is obtained that can be solved numerically.

New algorithms can in principle be designed by making the algorithmic parameters (such as $\alpha_k$  and $\beta_k$ in \cref{eq:nest1})  variable as well, and choosing these variables to minimize the value of the  maximization problem described above. 
In this way, \cite{KimFessler2016} derived new accelerated gradient by solving this min-max problem analytically.
Taylor, Hendrickx, and Glineur \cite{TaylorHendrickxGlineur2017} made significant enhancements to the PEP approach, including by developing precise finite parametrizations of the smooth convex and strongly convex function classes, tightening the PEP formulation, and studying the optimality of parametrization in several methods.

A related approach to analysis and construction of accelerated algorithms was proposed in \cite{LessardRechtPackard2016}.
They  view the iteration formulae as a dynamical system and deploy techniques from stability analysis of such systems.
Verification of convergence rates is performed by solving a semidefinite program (a matrix optimization problem in which the matrix variables are constrained to be symmetric and positive semidefinite). 
They apply the technique to steepest descent, heavy-ball, and Nesterov accelerated gradient.
In particular, they show that Polyak's choice of heavy-ball parameters $\alpha_k$ and $\beta_k$ in \eqref{eq:hb}, despite being optimal for strongly convex quadratics, do not even necessarily yield convergence to the unique solution for a strongly convex nonlinear function, for $\kappa>18$.
An extension of this work \cite{HuLessard2017} uses  {\em dissipativity}, another tool from control theory, as  a means to construct Lyapunov functions and thus derive convergence rates. 
With the algorithmic mechanism written as a dynamical system, quadratic constraints are introduced on the states and inputs in this system, based on known properties of the function (such as strongly convexity, Lipschitz smoothness of the gradient, co-coercivity) and the algorithm. 
A semidefinite program is then composed from  this information and solved to find a matrix $P$ that defines the Lyapunov function; see, for example \cite[Theorem~3]{HuWL18}.
This semidefinite program is essentially the dual of those encountered in the PEP approach.
The dissipativity approach can be deployed to analyze rates for general convex functions \cite[Section~4]{HuLessard2017} and for stochastic gradient methods \cite{HuWL18}, in which only an unbiased estimate of the gradient $\nabla f$ is available.

These approaches and the relationships between them are discussed in \cite[Appendix~C]{dAspremontScieurTaylor2021}.

\paragraph{Accelerated Gradient Methods for Nonconvex $f$.}
For gradient descent \eqref{eq:gdc}, we demonstrated above, using elementary analysis, the $O(1/\sqrt{T})$ rate for the quantity $\min_{k=0,1,\dotsc,T-1} \, \| \nabla f(x^k) \|$ for smooth nonconvex $f$ \cref{eq:grad3}. 
Although faster rates are not available in general, \cite[Section~6.4.1]{Lan20book} describes an accelerated gradient method that achieves the same $O(1/\sqrt{T})$ if $f$ is nonconvex, but has the accelerated $f(x^T)-f(x^*) = O(1/T^2)$ rate for convex $f$\cite[Corollary~6.14]{Lan20book}.
it requires knowledge of the Lipschitz constant for $\nabla f$ but no prior knowledge about whether or not $f$ is convex.

For finding convergence to approximate second-order points (satisfying \cref{eq:approx2o}), we already discussed the approach of \cite{JinNJ17a}, which is based on an accelerated gradient method.
Accelerated gradients are used in a different way to achieve the same goal by  \cite{carmon2018accelerated}. The rather complicated approach in this paper combines two strategies, one which uses the Lanczos method applied to the Hessian $\nabla^2 f(x^k)$ to identify directions of negative curvature for the Hessian, and one which applies accelerated gradient to a modified version of $f$.
Like \cite{JinNJ17a}, it has an iteration complexity of $O(\epsilon^{-7/4})$, a complexity also shared by the Newton-based methods discussed below. 
Like the latter, it requires calculation of Hessian-vector products ($\nabla^2 f(x)v$ for given $v$), so in this sense is not simply a first-order method.

\paragraph{Extensions.} 
Accelerated gradient methods can be adapted and extended to more general situations than the unconstrained problem \cref{eq:f} that we focus on here.
When the Lipschitz parameter $L$ is unknown, it can be estimated easily by a process similar to backtracking line search.
The methods can be adapted to constrained optimization problems of the form $\min_{x \in \Omega} f(x)$, for closed convex $\Omega$, provided that the projection operation $P_{\Omega}(x) := \arg\min_{y \in \Omega} \, \|x-y \|$ can be computed efficiently.
They can be extended too to regularized objectives such as $f(x) + \lambda \|x\|_1$ for some parameter $\lambda>0$, most famously in the FISTA method of Beck and Teboulle \cite{BecT08a}.
More generally, functions of the form $\min \, f(x) + h(x)$ can be minimized by accelerated methods, provided that $h(x)$ is a closed convex function whose prox-operator\footnote{The prox-operator is defined by $\mbox{prox}_h(v) := \arg\min_x \left\{ h(x) + \tfrac12 \|x-v\|_2^2 \right\}$.}  can be computed efficiently \cite{ParB13a}.
In these cases of constraints and regularized objectives, the straightforward extensions available for accelerated gradient methods lie in contrast to nonlinear conjugate gradient, which are difficult to extend beyond the unconstrained case.

There are accelerated-gradient versions of the mirror descent algorithm  \cite[Chapter~3]{NemY83},  used for minimization over the simplex $\Omega = \{ x \in \R^n \, : \, x \ge 0, \; \bfone^Tx =1 \}$ and applied to problems involving probability spaces. 

A technique that can improve convergence of gradient methods, both in theory and in practice, is {\em restarting}, which simply means that we periodically re-initialize the algorithm  at the  current iterate and proceed.
(We mentioned already that restarting improved the performance of first-order methods for linear programming.)
Restarting schemes can be devised for functions satisfying {\em Holderian error bounds}: $\frac{\mu}{r} \mbox{dist}(x,\cS)^r \le f(x)-f(x^*)$ for some positive $\mu$ and $r$, where $\cS$ is the set of minimizers of $f$. 
(Strongly convex functions \cref{eq:sconvex} satisfy these bounds with $r=2$.) 
For such functions, judicious restarting of schemes for general convex functions can yield convergence rates faster than the $O(1/T^2)$ rate of optimal gradient methods.
In fact, in the strongly convex case, rates matching those of accelerated gradient methods ($f(x^T)-f(x^*) \le \epsilon$ in $O(\sqrt{\kappa} | \log \epsilon |)$ iterations) can be attained \cite[Theorem~6.2]{dAspremontScieurTaylor2021}.
See \cite[Chapter~6]{dAspremontScieurTaylor2021} for further details of restarting schemes.

\subsection{Theory and Practice of Gradient Methods.}

Although far from exhaustive, our survey of  gradient-based methods for \eqref{eq:f} has included many techniques of elaborate design and a wealth of underlying theory. 
What are the relationships between this body of theory and current practice in nonlinear optimization?
Some algorithms, such as FISTA for $\ell_1$-regularized optimization \cite{BecT08a}, work well ``out of the box,'' matching their solid theoretical design with good practical performance on compressed sensing and related problems.
Nonlinear conjugate gradient methods were motivated by fascinating theory for the case of quadratic $f$, and some variants work well in practice on nonlinear $f$, but gaps remains with the weaker convergence theory for the latter case.
Restarting is an area in which theoretical and practical developments appear to have gone hand-in-hand, with computational experience being gathered as results were being proved.
The fascinating theory of versions of gradient descent \cref{eq:gd} that optimize progress over multiple steps, such as \cite{altschuler2025acceleration} does not have immediate practical relevance, though it provides insights into the possibilities of nonmonotone algorithms that take a longer view than the one-step perspective usually associated with gradient descent.

The fundamental idea of momentum that underlies the accelerated gradient methods described above has been enormously influential, though the great number of specific algorithms based on this idea are rarely implemented as written.
Momentum translates into each step being a weighted average of all gradients encountered so far (or since the latest restart), so the philosophy underlying momemtum and gradient averaging is the same.

The idea of momentum carries over to stochastic gradient methods for finite-sum objectives \cref{eq:finite-sum}, in which the gradient $\nabla f(x^k)$ is approximated by
\begin{equation} \label{eq:minib}
\frac{1}{|\cB_k|} \sum_{i \in \cB_k} \nabla f_i(x^k),
\end{equation}
where $\cB_k \subset \{1,2,\dotsc,N \}$ is a ``minibatch'' with $| \cB_k| \ll N$, chosen at random, so that \cref{eq:minib} is an {\em unbiased} estimate of $\nabla f(x^k)$.
Classical stochastic-gradient (``SGD'') methods \cite{RobM51} take steps along directions \cref{eq:minib} at each iteration, with steplengths chosen to guarantee convergence in expectation (see for example \cite{NemJLS09a} or \cite[Chapter~5]{WrightRecht22}).
The currently dominant variant of SGD is Adam \cite{Adam_2014}.
Adam adds two devices to the basic SGD scheme: weighted averaging over gradients encountered at previous iterations, and elementwise scaling of the components of the search direction, with the scaling also being learned over all previous iterations.
Given the enormous reach of Adam in the training of neural networks, it can be said that the idea of momentum in gradient methods is  having a great impact on the practice of optimization in machine learning applications.


\subsection{Newton's Method.} \label{sec:newton}

Newton's method for \cref{eq:f}, in its essential form, generates iterates as follows (starting from a given point $x^0$):
\begin{quote}
     Solve $\nabla^2 f(x^k) s^k = -\nabla f(x^k)$, \ $x^{k+1} \leftarrow x^k + s^k$,  \ $k=0,1,2,\dotsc.$
\end{quote}
For strictly convex $f$ (for which $\nabla^2 f(x^k)$ is positive definite for all $k$),  each step $s^k$ is the minimizer of the quadratic Taylor-series expansion of $f$ around $x^k$, which is
\begin{equation} \label{eq:newt1}
  m(s;x^k) := f(x^k) + \nabla f(x^k)^T s + \frac12 s^T \nabla^2 f(x^k) s.
\end{equation}
Various modifications can be applied to this basic approach to make it suitable for smooth nonlinear problems, even nonconvex problems. 
For example, line searches or trust regions can be applied at each step to ensure convergence to (or accumulation at) points satisfying optimality conditions, or the {\em Newton equations} $\nabla^2 f(x^k) s^k = -\nabla f(x^k)$ can be solved inexactly, possibly by an iterative scheme like linear CG, that requires only computation of matrix-vector products of the form $\nabla^2 f(x^k) p$ for some $p \in \R^n$.
Such products may not require  explicit computation and storage of the Hessian $\nabla^2 f(x^k)$, as we discuss below.

The best-known property of Newton's method is its asymptotic quadratic convergence. 
That is, once the iterates reach a neighborhood of a point $x^*$ satisfying the second-order sufficient conditions $\nabla f(x^*)=0$, $\nabla^2 f(x^*) \succ 0$, we have $\|x^{k+1}-x^* \| \le C \| x^k-x^*\|^2$, for some $C>0$, for the remainder of the sequence.
Because this convergence is so rapid, the theory is relevant only on the last few iterates of a typical run of Newton's method. 
But for some kinds of problems, the extra expense associated with using second derivatives is more than compensated for by faster convergence of the iterates, even far from solutions.

For smooth nonconvex $f$, the nonasymptotic rate of Newton's method is no faster than gradient descent.
\cite[Section~3.1]{cartis2022evaluation} constructs a smooth function $f:\R^2 \to \R$ for which, given an accuracy threshold $\epsilon_1>0$, Newton's method may require more than $1/\epsilon_1^{2}$ iterations to identify a point $x$ with $\| \nabla f(x) \| \le \epsilon_1$.
Although this example is compelling, the construction of $f$ depends on the choice of $\epsilon_1$; it is unclear whether there exists a function $f$ for which the basic Newton's method above requires more than $1/\epsilon^{2}$ iterations for an arbitrary $\epsilon$.
Indeed, progress over the past 20 years has yielded Newton-based methods for which the complexity is considerably better than this example would suggest. 
These methods work for nonconvex functions, and can find points that approximately satisfy not just the first-order condition $\nabla f(x)=0$, but also the second-order sufficient condition $\nabla^2 f(x) \succeq 0$.

A trust-region variant of Newton's method, described in \cite[Chapter~4]{NocW06}, solves at each iteration a subproblem of the form
\begin{equation} \label{eq:newt2}
s^k :=   \arg\min_{s \in \R^n} \, m(s;x^k) \;\; \mbox{subject to} \; \| s\|_2^2 \le \Delta_k,
\end{equation}
for $m(s;x^k)$ defined in \cref{eq:newt1} and some {\em trust-region radius} $\Delta_k$, which is increased and decreased according to whether the step $s^k$ from \cref{eq:newt2} yields a sufficiently large decrease in the objective function $f$.
This subproblem can be solved efficiently, even when the Hessian $\nabla^2 f(x^k)$ is indefinite.
Asymptotic analysis of different variants of the trust-region approach can be shown to generate sequences $\{x^k \}$ whose accumulation points are stationary, or even satisfy second-order conditions.
In \cite{CurRS14a},  a  trust-region variant (motivated by the cubic regularization method discussed below) was shown to require $O(\epsilon^{-3/2})$ iterations to  converge to a point with $\| \nabla f(x) \| \le \epsilon$.
A later approach, hewing more closely to the traditional trust-region method but with modifications that equip it with nonasymptotic complexity properties, is described in \cite{Cur19a}.

The {\em cubic regularization} method \cite{NesP06a} was possibly the
first Newton-based method devised specifically to admit nonasymptotic
complexity theory, while also having practical relevance because of its relationship to the trust-region approach.
This method assumes that the Hessian $\nabla^2 f(x)$ is Lipschitz continuous with Lipschitz parameter $M$ \cref{eq:hlip}.
 $f$ can thus be bounded by a cubic function, as follows:
\begin{equation} \label{eq:cubic}
  f(x^k+s)  \le t(s;x^k,M) := m(s;x^k) + \frac{M}{6} \|s\|^3.
\end{equation}
The cubic regularization algorithm takes the step $s^k$ to be the
minimizer of $t(s;x^k,M)$, and sets $x^{k+1} = x^k+s^k$.
Within $T = O(\epsilon^{-3/2})$ iterations, at least one
of the iterates $k=1,2,\dotsc,T$ will identify a point that is
$\epsilon$-approximately second-order necessary, in the sense  of \cref{eq:approx2o}, 
where $\epsg = O(\epsilon)$ and $\epsH = O(\sqrt{\epsilon})$.
In fact, this property holds for a more adaptive version of the algorithm, explicated in \cite{CarGT11a,CarGT11b}, in which $M$ is replaced by a parameter that ensures sufficient reduction of the function $f$.
(A method like this had been proposed much earlier by Griewank \cite{Gri81b}, without nonasymptotic analysis.)
Finding the minimizer of $t(s;x^k,M)$ is not straightforward, but there is a close relationship between this minimizer and the solution of the trust-region problem \cref{eq:newt2} (as can be shown by an elementary application of the Karush-Kuhn-Tucker  theory for constrained optimization).

Practical Newton methods  for large-scale problems make use of {\em inexact} solutions of the Newton equations $\nabla^2 f(x^k) s^k = -\nabla f(x^k)$, often by applying the conjugate gradient method for minimizing convex quadratics to the quadratic approximation \eqref{eq:newt1}.
When the Hessian $\nabla^2 f(x^k)$ is not positive definite, the conjugate gradient method is usually quick to identify a direction $d^k$ of negative curvature, that is, $(d^k)^T \nabla^2 f(x^k)d^k<0$.
Such a direction is still useful as a search direction that reduces the function significantly.
A trust-region algorithm that exploits this and other properties of the conjugate gradient methods was described in \cite{Ste83}; see also \cite[Section~7.1]{NocW06}.
The method described in \cite{Cur19a} is a modification of this approach that equips it with non-asymptotic properties, including convergence to a point satisfying conditions \cref{eq:approx2o} in $\tilde{O}(\epsilon^{-3/2})$ iterations.
``Line-search'' versions of the inexact Newton method are also available, depending on their asymptotic global convergence properties not a trust region of the form $\|s\| \le \Delta$, but rather on the choice of a steplength $\alpha$ along the search direction $s^k$ (or $d^k$)  that ensures a significant decrease in the objective function at each step.
A modification of the line-search approach that equips it with non-asymptotic properties is described in \cite{RoyOW18a}.

The iteration complexity results of these Newton-based algorithms can be converted to oracle complexities.
If we define the ``oracle'' to include explicit evaluation of $\nabla^2 f(x)$, as well as the gradient $\nabla f(x)$ and function value $f(x)$, the transition is trivial.
A more challenging task is to devise oracle complexities that rely only on evaluation of first-order derivatives of $f$.
The key observation here is that the Hessian-vector product $\nabla^2 f(x)p$ for any $x, p \in \R^n$ can be obtained by differentiating the inner product $\nabla f(x)^Tp$ with respect to $x$.
If reverse-mode automatic differentiation is used, the cost is bounded by five times the cost of evaluating the gradient $\nabla f(x)$ itself; see for example \cite{GriW08}.
Thus, a kind of oracle complexity can be derived by counting the number of gradient evaluations plus the number of Hessian-vector products, and adding the cost of evaluating the function $f$ (usually cheaper). 
The algorithms in \cite{RoyOW18a,Cur19a} require $\tilde{O}(\epsilon^{-7/4})$ such operations in total.
An important element of the proof of this result is to show that no more than $O(\epsilon^{-1/2})$ iterations of conjugate gradient are required within each iteration of the modified Newton method in order to compute a direction that leads to significant progress in reducing $f$.

Several more things are worth mentioning about the methods of
\cite{RoyOW18a,Cur19a}.  They are based on practical methods, but the
modifications that are made to admit nonasymptotic theory do not improve the
practical performance. Moreover, the complexity bounds are quite
pessimistic for small $\epsilon$; there remains a large gap between
these bounds and practical performance. This is not surprising since
our assumptions on $f$ are so mild, and potentially admit pathological
examples.
Finally, we note that the ``oracle'' complexity bound of $\tilde{O}(\epsilon^{-7/4})$ achieved by these methods matches that of the perturbed accelerated gradient method in \cite{JinNJ17a}, which does not make use of Hessians at all.

\paragraph{Newton's Method on Self-Concordant Barrier Functions.}
Newton's method is a key component of the primal path-following methods for conic programming described in \cite[Chapter~3]{NesN93}.
These methods apply to constrained optimization problems of the form
\begin{equation} \label{eq:fomega}
 \min_{x \in \R^n} \, c^Tx \quad \mbox{subject to} \; Ax=b, \quad x \in \Omega,
\end{equation}
where $\Omega \subset \R^n$ is a closed cone admitting a self-concordant barrier function.
The special properties of the  self-concordant barrier make the behavior of the Newton step more predictable over a wider region than is the case of general functions $f$. 
Newton's method, either with unit steps or with steplengths calculated by a prescribed formula (see \cite[p.~77]{NesN93}), can be applied in a systematic way to generate a sequence of points that follow the central path to a solution of \cref{eq:fomega}.
These methods find $\epsilon$-accurate solutions with iteration complexity $O(\vartheta |\log \epsilon|)$ (long-step variant) or $O(\sqrt\vartheta |\log \epsilon|)$ (short-step variant), where $\vartheta$ is the parameter of the barrier function (see \cite[Definition~2.3.1]{NesN93}).
For the nonnnegative orthant $\{ x \in \R^n \, : \, x \ge 0 \}$, $-\sum_{i=1}^n \ln x_i$ is a self-concordant barrier function with parameter $n$.


\subsection{Quasi-Newton Methods.}

Starting in 1959 \cite{davidon1959variable}, methods have been proposed that are inspired by Newton's method, but which form and maintain approximations of the Hessian $\nabla^2 f(x)$ (or its inverse) from gradients $\nabla f(x^k)$ encountered during the algorithm. 
Taylor's theorem provides the key observation: Given any successive pair of iterates $x^{k}$, $x^{k+1}$ with $s^k := x^{k+1}-x^{k}$ and $y^k := \nabla f(x^{k+1})-\nabla f(x^{k})$, we have that $\nabla^2 f(x^{k+1}) s^k \approx y^k$.
Quasi-Newton methods define matrices $B_{k+1}$ that capture this property of Hessian.
In fact $B_{k+1}$ is required to satisfy the condition $B_{k+1} s^k = y^k$ (known as the {\em secant equation}), along with other desirable properties, such as positive definiteness, and changing minimally in some sense from the previous matrix $B_k$. 
The quasi-Newton search direction at $x^{k+1}$ is then obtained by solving $B_{k+1} s^{k+1} = -\nabla f(x^{k+1})$, and defining the new iterate to be $x^{k+2} = x^{k+1} + \alpha_{k+1} s^{k+1}$.

An alternative approach, easier to implement, maintains an approximation $H_{k+1}$ to the {\em inverse} Hessian instead, and requires it to satisfy the appropriate form of the secant equation: $s^{k+1} = H_{k+1} y^{k+1}$, along  with other desirable properties. The step $s^{k+1}$ can  then be calculated by a matrix multiplication: $s^{k+1} = -H_{k+1} \nabla f(x^{k+1})$.

The most popular quasi-Newton method is the BFGS method, proposed independently in 1970 by four inventors (Broyden, Fletcher, Goldfarb, and Shanno). 
Provided that $(s^k)^T y^k >0$ (a property that can be guaranteed if approximate optimality conditions are imposed on the steplength $\alpha_k$), there is an explicit formula for updating $H_k$ to obtain $H_{k+1}$:
\begin{equation}
    \label{eq:bfgs}
    H_{k+1} = (I-\rho_k s^k (y^k)^T) H_k (I - \rho_k y^k (s^k)^T) + \rho_k s^k (s^k)^T, \quad \mbox{where $\rho_k := 1/(s^k)^T y^k$.}
\end{equation}
The method is often initialized by a scalar multiple of the identity, that captures the appropriate scaling of the  true Hessian inverse: $H_0 = \left[ (y^0)^Ts^0/(y^0)^T y^0 \right] \,  I$. 

The best-known property of quasi-Newton methods is their superlinear local convergence, proved for inexact line searches in the 1970s  \cite{broyden1973local,dennis1974characterization}, and observed in practice.
Asymptotic global convergence for BFGS with inexact line searches on convex functions was proved by Powell in 1976 \cite{Powell76}. (For nonconvex functions, \cite{Dai02} showed that the algorithm may cycle.)
Non-asymptotic convergence results have taken much longer to emerge.
A bound of $O(\kappa | \log \epsilon |)$ for  BFGS with inexact line searches on convex functions is proved in \cite[Theorem~4.1]{jin2025nonasymptoticglobalconvergenceanalysis}, matching the complexity of steepest descent, but with a constant factor that is potentially large (and dependent on the parameters in the inexact line search conditions). 
A second result in \cite{jin2025nonasymptoticglobalconvergenceanalysis} shows a linear rate that is independent of the condition number, but this rate cuts in only after $O(\kappa)$ iterations.

Although they use only first-order information, quasi-Newton methods typically require much more storage ($O(n^2)$, to store the matrices $H_k$ or $B_k$) than the gradient, conjugate gradient, and accelerated gradient methods discussed above, which require storage of just a few vectors of length $n$.
The limited-memory BFGS method of \cite{LiuN89} (known universally as ``L-BFGS'') reduces these storage requirements dramatically by maintaining update vectors $(s^k,y^k)$ from only the past $m$ iterations, where $m \ll n$. (Values $m=5$ or $m=10$ are typical.) 
At iteration $j+1$, initializing with $H_{j+1}^0 = \left[ (y^j)^Ts^j/(y^j)^T y^j \right] \,  I$, the update formulae \cref{eq:bfgs} is applied for $k=j-m+1, j-m+2, \dotsc, j$ to obtain the approximation $H_{j+1}$, with the new search direction being set to $s^{j+1} = -H_{j+1} \nabla f(x^{j+1})$. 
Importantly, the matrix $H_{j+1}$ does not need to be stored explicitly; $s^{j+1}$ can be calculated efficiently in $O(mn)$ operations from the vectors $(s^k,y^k)$ for $k=j-m+1, \dotsc, j$.
Theory for L-BFGS is limited. 
A safeguarded version can be shown to have similar asymptotic convergence properties as gradient descent, and superlinear local convergence is not seen in general.
However, L-BFGS is widely used and effective in practice, often outperforming nonlinear conjugate gradient.

We note that one variant of  the Barzilai-Borwein method \cite{BarB88} can be viewed as simply L-BFGS with no updates ($m=0$) and no line search. 
The initialization $H_{j+1}^0 = \left[(y^j)^Ts^j/(y^j)^T y^j \right] \, I$ used in L-BFGS is the actual Hessian approximation used by Barzilai-Borwein.

\subsection{``Tractable'' Nonconvex Problems.} \label{sec:nonconvex}

The problem of finding the global minimizer of a general smooth function $f:\R^n \to \R$ is intractable in general.
Most algorithms for such functions (including those described above) can be proved only to converge to points that approximately satisfy the first-order optimality conditions or sometimes second-order conditions (see \cref{eq:approx2o}), which are markers only of {\em local} solutions at best.

In recent years, however, many problem classes have been identified (largely from the areas of matrix optimization and machine learning) for which finding global optima is a tractable problem.
(See \cite{Sun21} for a list of such problem classes.)
In some cases (such as those involving nonconvex formulations of low-rank matrix optimization problems) this phenomenon was observed first in practice: Algorithms that were guaranteed to find only first-order optimal points were consistently finding global solutions.
Explanations of why this happens, including rigorous specification of conditions under which it happens, has followed.

We discuss just the first and most famous instance of this phenomenon, in semidefinite programming, explicated by Burer and Monteiro \cite{BurM03a,BurM05a}.
Semidefinite programming is the canonical matrix optimization problem, whose standard form is
\begin{equation}
    \label{eq:sdp}
    \min_{X \in S\R^{n \times n}} \, \langle C,X \rangle \quad \mbox{subject to $\langle A_i,X \rangle = b_i$, \ $i=1,2,\dotsc,m$; \ \  $X \succeq 0$;}
\end{equation}
where the variable $X$ and the data $C$ and $A_i$, $i=1,2,\dotsc,m$ are symmetric $n \times n$ matrices, and the inner product is defined by $\langle X,Y \rangle = \sum_{i,j=1}^n X_{ij} Y_{ij}$.
This is a convex conic optimization problem, and interior-point methods with polynomial complexity have been well known since the 1990s.
However, these approaches are computationally impractical for large $n$.

There are solutions  $X$ to this problem whose rank is bounded approximately by $\sqrt{2m}$, so by choosing $r$ greater than this threshold,  \cref{eq:sdp} can be reformulated by writing $X=ZZ^T$, where $Z \in \R^{n \times r}$ is a tall matrix with $r \ll n$. 
This factored form of $X$ has  positive semidefiniteness ``built in,'' so we can omit the constraint $X \succeq 0$ and substitute in \cref{eq:sdp} to obtain
\begin{equation}
    \label{eq:sdp.low}
    \min_{Z \in \R^{n \times r}} \, \langle C,ZZ^T \rangle \quad \mbox{subject to $\langle A_i,ZZ^T \rangle = b_i$, \ $i=1,2,\dotsc,m$.}
\end{equation}
This formulation is nonconvex, and  can be solved using algorithms for nonlinear nonconvex  constrained optimization.
An augmented Lagrangian method with regularization for \cref{eq:sdp.low} is described in \cite{BurM05a} and the sequence of iterates $Z^k$ generated by this approach is shown to have the property that the accumulation points of $\{ Z^k (Z^k)^T \}$ are optimal for \cref{eq:sdp}.

The Burer-Monteiro parametrization and nonconvex formulation has found widespread use because of its practical efficiency resulting from the vastly reduced dimensionality. 
The theory underlying this approach has continued to develop since the initial papers. A recent result \cite{CifM22} establishes polynomial-time guarantees in the smoothed analysis setting, in which a random perturbation is added to the cost matrix $C$.

\section{Discussion.} \label{sec:discussion}

We have discussed  theoretical results for algorithms in two areas of optimization: linear programming and smooth convex optimization.
We have also indicated which algorithms and approaches are most widely used in practice, and discussed how theory and practice have influenced each other in the development of algorithms.
In both areas, innovations continue apace in both theory and practice. 
To mention one example: Interior-point methods for LP have become a focus of intense interest in the combinatorial optimization community in recent years.

We have left much ground uncovered: Other areas of optimization have seen enormous developments in theory and practice in recent years, with similar interplay between these two aspects of the discipline. 
We mention in particular semidefinite programming, constrained  nonlinear optimization, and stochastic-gradient (SGD) methods, which are  based on gradient approximations \cref{eq:minib} for finite-sum problems \cref{eq:finite-sum}.
The latter topic warrants a full article of its own! 
Its current importance cannot be overstated: SGD is undoubtedly consuming far more compute cycles than any other optimization algorithm in the world today, via its importance in training neural networks, including LLMs. 
Among {\em all} computing on CPUs and GPUs, it likely ranks second in total cycles used, exceeded only by neural network inference.
It is possibly surprising that an algorithm based on crude (but cheap) approximations to the gradient admits any useful theory at all --- but it does.
Dating to 1951 \cite{RobM51}, the basic theory of the method is outlined in \cite{NemJLS09a} (with many relevant works published before and since). 
Recent developments of the basic algorithm have focused on the use of momentum / gradient averaging and diagonal scaling \cite{DucHS10a,Adam_2014}, the choice of steplength schedules, and parallel implementation (for example \cite{BerT89,Hogwild} and many works on federated learning).
Theoretical developments have taken place alongside a huge amount of computational experimentation, with influences flowing in both directions. 
(Theoreticians strive to justify approaches that have seen practical success; practitioners gain confidence from theoretical underpinnings of  innovations that have proved useful in practice.)








\section*{Acknowledgments.}

Thanks to Ben Recht for suggesting the topic of this review and to Haihao Lu for his guidance on first-order methods for linear programming.

\bibliographystyle{siamplain}
\bibliography{icm-wright,refs}

\begin{thebibliography}{100}

\bibitem{alman2024asymmetryyieldsfastermatrix}
{\sc J.~Alman, R.~Duan, V.~V. Williams, Y.~Xu, Z.~Xu, and R.~Zhou}, {\em More
  asymmetry yields faster matrix multiplication}, in Proceedings of the 2025
  Annual ACM-SIAM Symposium on Discrete Algorithms (SODA), SIAM, 2025,
  pp.~2005--2039.

\bibitem{altschuler2018greed}
{\sc J.~M. Altschuler}, {\em Greed, hedging, and acceleration in convex
  optimization}, master's thesis, Massachusetts Institute of Technology,
  Cambridge, MA, 2018.

\bibitem{altschuler2023acceleration}
{\sc J.~M. Altschuler and P.~A. Parrilo}, {\em Acceleration by stepsize hedging
  {II}: Silver stepsize schedule for smooth convex optimization}, arXiv
  preprint arXiv:2309.16530,  (2023).
\newblock To appear in {\em Mathematical Programming, Series A}.

\bibitem{altschuler2025acceleration}
{\sc J.~M. Altschuler and P.~A. Parrilo}, {\em Acceleration by stepsize
  hedging: Multi-step descent and the silver stepsize schedule}, Journal of the
  ACM, 72 (2025), pp.~1--38.

\bibitem{applegate2021practical}
{\sc D.~Applegate, M.~D{\'\i}az, O.~Hinder, H.~Lu, M.~Lubin, B.~O'Donoghue, and
  W.~Schudy}, {\em Practical large-scale linear programming using primal-dual
  hybrid gradient}, in Proceedings of the 35th International Conference on
  Neural Information Processing Systems, 2021, pp.~20243--20257.

\bibitem{applegate2025pdlp}
{\sc D.~Applegate, M.~D{\'\i}az, O.~Hinder, H.~Lu, M.~Lubin, B.~O'Donoghue, and
  W.~Schudy}, {\em {PDLP:} a practical first-order method for large-scale
  linear programming}, arXiv preprint arXiv:2501.07018,  (2025).

\bibitem{bach2025optimalsmoothedanalysissimplex}
{\sc E.~Bach and S.~Huiberts}, {\em Optimal smoothed analysis of the simplex
  method}, in Proceedings of the 66th IEEE Symposium on Foundations of Computer
  Science, 2025, \url{https://arxiv.org/abs/2504.04197},
  \url{https://arxiv.org/abs/2504.04197}.

\bibitem{BarB88}
{\sc J.~Barzilai and J.~M. Borwein}, {\em Two-point step size gradient
  methods}, IMA Journal of Numerical Analysis, 8 (1988), pp.~141--148.

\bibitem{BecT08a}
{\sc A.~Beck and M.~Teboulle}, {\em A fast iterative shrinkage-threshold
  algorithm for linear inverse problems}, SIAM Journal on Imaging Sciences, 2
  (2009), pp.~183--202.

\bibitem{BerT89}
{\sc D.~P. Bertsekas and J.~N. Tsitsiklis}, {\em Parallel and Distributed
  Computation: Numerical Methods}, Prentice-Hall, Inc., Englewood Cliffs, New
  Jersey, 1989.

\bibitem{blum2012complexity}
{\sc L.~Blum, F.~Cucker, M.~Shub, and S.~Smale}, {\em Complexity and Real
  Computation}, Springer Science \& Business Media, 2012.

\bibitem{Bor87}
{\sc K.-H. Borgwardt}, {\em The Simplex Method, A Probabilistic Analysis},
  Springer-Verlag, Berlin, 1987.

\bibitem{Borgwardt1999}
{\sc K.-H. Borgwardt}, {\em A sharp upper bound for the expected number of
  shadow vertices in {LP}-polyhedra under orthogonal projection on
  two-dimensional planes}, Mathematics of Operations Research, 24 (1999),
  pp.~544--603, \url{https://doi.org/10.1287/moor.24.3.544}.

\bibitem{broyden1973local}
{\sc C.~G. Broyden, J.~E. Dennis, and J.~J. Mor{\'e}}, {\em On the local and
  superlinear convergence of quasi-newton methods}, Journal of the Institute of
  Mathematics and its Applications, 12 (1973), pp.~223--245.

\bibitem{BurM03a}
{\sc S.~Burer and R.~D.~C. Monteiro}, {\em A nonlinear programming algorithm
  for solving semidefinite programs via low-rank factorizations}, Mathematical
  Programming, Series {B}, 95 (2003), pp.~329--257.

\bibitem{BurM05a}
{\sc S.~Burer and R.~D.~C. Monteiro}, {\em Local minima and convergence in
  low-rank semidefinite programming}, Mathematical Programming, Series {A}, 103
  (2005), pp.~427--444.

\bibitem{carmon2018accelerated}
{\sc Y.~Carmon, J.~C. Duchi, O.~Hinder, and A.~Sidford}, {\em Accelerated
  methods for nonconvex optimization}, SIAM Journal on Optimization, 28 (2018),
  pp.~1751--1772.

\bibitem{cartis2022evaluation}
{\sc C.~Cartis, N.~I. Gould, and P.~L. Toint}, {\em Evaluation Complexity of
  Algorithms for Nonconvex Optimization: Theory, Computation and Perspectives},
  SIAM, 2022.

\bibitem{CarGT11a}
{\sc C.~Cartis, N.~I.~M. Gould, and P.~L. Toint}, {\em Adaptive cubic
  regularisation methods for unconstrained optimization. {Part I:} motivation,
  convergence and numerical results}, Mathematical Programming, Series {A}, 127
  (2011), pp.~245--295.

\bibitem{CarGT11b}
{\sc C.~Cartis, N.~I.~M. Gould, and P.~L. Toint}, {\em Adaptive cubic
  regularisation methods for unconstrained optimization. part ii: worst-case
  function-and derivative-evaluation complexity}, Mathematical Programming,
  Series {A}, 130 (2011), pp.~295--319.

\bibitem{chambolle2011first}
{\sc A.~Chambolle and T.~Pock}, {\em A first-order primal-dual algorithm for
  convex problems with applications to imaging}, Journal of mathematical
  imaging and vision, 40 (2011), pp.~120--145.

\bibitem{CifM22}
{\sc D.~Cifuentes and A.~Moitra}, {\em Polynomial time guarantees for the
  burer-monteiro method}, in Advances in Neural Information Processing Systems,
  vol.~35, 2022, pp.~23923--23935.

\bibitem{cohen2021solving}
{\sc M.~B. Cohen, Y.~T. Lee, and Z.~Song}, {\em Solving linear programs in the
  current matrix multiplication time}, {Journal of the ACM (JACM)}, 68 (2021),
  pp.~1--39.

\bibitem{Cur19a}
{\sc F.~E. Curtis, D.~P. Robinson, C.~W. Royer, and S.~J. Wright}, {\em
  Trust-region {Newton-CG} with strong second-order complexity guarantees for
  nonconvex optimization}, {SIAM} Journal on Optimization, 31 (2021),
  pp.~518--544, \url{https://arxiv.org/abs/1912.04365}.

\bibitem{CurRS14a}
{\sc F.~E. Curtis, D.~P. Robinson, and M.~Samadi}, {\em A trust-region
  algorithm with a worst-case iteration complexity of ${O}(\epsilon^{-3/2})$
  for nonconvex optimization}, Mathematical Programming, Series {A}, 162
  (2017), pp.~1--32.

\bibitem{Dai02}
{\sc Y.-H. Dai}, {\em Convergence properties of the bfgs algoritm}, SIAM
  Journal on Optimization, 13 (2002), pp.~693--701.

\bibitem{dant:histo}
{\sc G.~Dantzig}, {\em Linear programming, the story about how it began}, in
  History of Mathematical Programming: {A} Collection of Personal
  Reminiscences, J.~K. Lenstra, A.~H. G.~R. Kan, and A.~Schrijver, eds.,
  Amsterdam, 1991, Elsevier.

\bibitem{das2024nonlinear}
{\sc S.~Das~Gupta, R.~M. Freund, X.~A. Sun, and A.~Taylor}, {\em Nonlinear
  conjugate gradient methods: worst-case convergence rates via
  computer-assisted analyses}, Mathematical Programming,  (2024), pp.~1--49.

\bibitem{dAspremontScieurTaylor2021}
{\sc A.~d'Aspremont, D.~Scieur, and A.~Taylor}, {\em Acceleration methods},
  Foundations and Trends in Optimization, 5 (2021), pp.~1--245,
  \url{https://doi.org/10.1561/2400000046}.

\bibitem{davidon1959variable}
{\sc W.~C. Davidon}, {\em Variable metric method for minimization}, tech.
  report, Argonne National Laboratory, Lemont, Ill., USA, 1959.

\bibitem{dennis1974characterization}
{\sc J.~E. Dennis and J.~J. Mor{\'e}}, {\em A characterization of superlinear
  convergence and its application to quasi-newton methods}, Mathematics of
  Computation, 28 (1974), pp.~549--560.

\bibitem{Dik67}
{\sc I.~I. Dikin}, {\em Iterative solution of problems of linear and quadratic
  programming}, Soviet Mathematics-Doklady, 8 (1967), pp.~674--675.

\bibitem{disser2023exponential}
{\sc Y.~Disser, O.~Friedmann, and A.~V. Hopp}, {\em An exponential lower bound
  for {Zadeh}’s pivot rule}, Mathematical Programming, 199 (2023),
  pp.~865--936.

\bibitem{Drori2017}
{\sc Y.~Drori}, {\em The exact information-based complexity of smooth convex
  minimization}, Journal of Complexity, 39 (2017), pp.~1--16,
  \url{https://doi.org/10.1016/j.jco.2016.11.002}.

\bibitem{DroriTeboulle2014}
{\sc Y.~Drori and M.~Teboulle}, {\em Performance of first-order methods for
  smooth convex minimization: A novel approach}, Mathematical Programming, 145
  (2014), pp.~451--482, \url{https://doi.org/10.1007/s10107-013-0653-0}.

\bibitem{DucHS10a}
{\sc J.~Duchi, E.~Hazan, and Y.~Singer}, {\em Adaptive subgradient methods for
  online learning and stochastic optimization.}, Journal of machine learning
  research, 12 (2011).

\bibitem{FleR64}
{\sc R.~Fletcher and C.~M. Reeves}, {\em Function minimization by conjugate
  gradients}, Computer Journal, 7 (1964), pp.~149--154.

\bibitem{gilbert1992global}
{\sc J.~C. Gilbert and J.~Nocedal}, {\em Global convergence properties of
  conjugate gradient methods for optimization}, SIAM Journal on optimization, 2
  (1992), pp.~21--42.

\bibitem{GilMSTW86}
{\sc P.~E. Gill, W.~Murray, M.~A. Saunders, J.~A. Tomlin, and M.~H. Wright},
  {\em On projected {N}ewton barrier methods for linear programming and an
  equivalence to {K}armarkar's projective method}, Mathematical Programming, 36
  (1986), pp.~183--209.

\bibitem{GolV83}
{\sc G.~H. Golub and C.~F. {Van Loan}}, {\em Matrix Computations}, The Johns
  Hopkins University Press, {B}altimore, second~ed., 1989.

\bibitem{goujaud2023provable}
{\sc B.~Goujaud, A.~Taylor, and A.~Dieuleveut}, {\em Provable non-accelerations
  of the heavy-ball method}, arXiv preprint arXiv:2307.11291,  (2023).

\bibitem{Gri81b}
{\sc A.~Griewank}, {\em The modification of {Newton}'s method for unconstrained
  optimization by bounding cubic terms}, {Technical Report} NA/12, DAMTP,
  Cambridge University, 1981.

\bibitem{GriW08}
{\sc A.~Griewank and A.~Walther}, {\em Evaluating Derivatives: Principles and
  Techniques of Algorithmic Differentiation}, Frontiers in Applied Mathematics,
  SIAM, Philadelphia, PA, second~ed., 2008.

\bibitem{grimmer2024provably}
{\sc B.~Grimmer}, {\em Provably faster gradient descent via long steps}, SIAM
  Journal on Optimization, 34 (2024), pp.~2588--2608.

\bibitem{HagZ05}
{\sc W.~W. Hager and H.~Zhang}, {\em A new conjugate gradient method with
  guaranteed descent and efficient line search}, {SIAM} Journal on
  Optimization, 16 (2005), pp.~170--192.

\bibitem{HagZ06b}
{\sc W.~W. Hager and H.~Zhang}, {\em A survey of nonlinear conjugate gradient
  methods}, Pacific Journal of Optimization, 2 (2006), pp.~35--58.

\bibitem{Hestenes}
{\sc M.~R. Hestenes and E.~Stiefel}, {\em Methods of conjugate gradients for
  solving linear systems}, Journal of Research of the National Bureau of
  Standards, 49 (1952), pp.~409--436.

\bibitem{HuLessard2017}
{\sc B.~Hu and L.~Lessard}, {\em Dissipativity theory for nesterov's
  accelerated method}, in Proceedings of the 34th International Conference on
  Machine Learning (ICML), vol.~70 of Proceedings of Machine Learning Research,
  PMLR, 2017, pp.~1549--1557.

\bibitem{HuWL18}
{\sc B.~Hu, S.~J. Wright, and L.~Lessard}, {\em Dissipativity theory for
  accelerating stochastic variance reduction: {A} unified analysis of {SVRG}
  and katyusha using semidefinite programs}, in Proceedings of the 35th
  International Conference on Machine Learning, {ICML} 2018, 2018,
  pp.~2043--2052.

\bibitem{jiang2021faster}
{\sc S.~Jiang, Z.~Song, O.~Weinstein, and H.~Zhang}, {\em A faster algorithm
  for solving general lps}, in Proceedings of the 53rd Annual ACM SIGACT
  Symposium on Theory of Computing, 2021, pp.~823--832.

\bibitem{jin2017escape}
{\sc C.~Jin, R.~Ge, P.~Netrapalli, S.~M. Kakade, and M.~I. Jordan}, {\em How to
  escape saddle points efficiently}, in International conference on machine
  learning, PMLR, 2017, pp.~1724--1732.

\bibitem{JinNJ17a}
{\sc C.~Jin, P.~Netrapalli, and M.~I. {Jordan}}, {\em Accelerated gradient
  descent escapes saddle points faster than gradient descent}, Journal of
  Machine Learning Research, 75 (2018), pp.~1--44.

\bibitem{jin2025nonasymptoticglobalconvergenceanalysis}
{\sc Q.~Jin, R.~Jiang, and A.~Mokhtari}, {\em Non-asymptotic global convergence
  analysis of {BFGS} with the {Armijo-Wolfe} line search}, in The Thirty-eighth
  Annual Conference on Neural Information Processing Systems, 2024.

\bibitem{KarNS16a}
{\sc H.~Karimi, J.~Nutini, and M.~Schmidt}, {\em Linear convegence of gradient
  and proximal-gradient methods under the {Polyak-{\L}ojasiewicz} condition},
  in European Conference on Machine Learning, 2016.

\bibitem{KarV21a}
{\sc S.~Karimi and S.~A. Vavasis}, {\em Nonlinear conjugate gradient for smooth
  convex functions}, Mathematical Programming Computation, 16 (2024),
  pp.~229--254.

\bibitem{Kar84}
{\sc N.~Karmarkar}, {\em A new polynomial-time algorithm for linear
  programming}, Combinatorica, 4 (1984), pp.~373--395.

\bibitem{Kha79}
{\sc L.~G. Khachiyan}, {\em A polynomial algorithm in linear programming},
  Soviet Mathematics Doklady, 20 (1979), pp.~191--194.

\bibitem{KimFessler2016}
{\sc D.~Kim and J.~A. Fessler}, {\em Optimized first-order methods for smooth
  convex minimization}, Mathematical Programming, 159 (2016), pp.~81--107,
  \url{https://doi.org/10.1007/s10107-015-0949-3}.

\bibitem{Adam_2014}
{\sc D.~P. Kingma and J.~Ba}, {\em Adam: {A} method for stochastic
  optimization}, in Proceedings of International Conference on Learning
  Representations, 2015.

\bibitem{KleM72}
{\sc V.~Klee and G.~J. Minty}, {\em How good is the simplex algorithm?}, in
  Inequalities, O.~Shisha, ed., Academic Press, New York, 1972, pp.~159--175.

\bibitem{Knu91}
{\sc D.~E. Knuth}, {\em Theory and practice}, Theoretical Computer Science, 90
  (1991), pp.~1--15.

\bibitem{Lan20book}
{\sc G.~Lan}, {\em First-order and stochastic optimization methods for machine
  learning}, Springer Series in the Data Sciences, Springer-Nature, 2020.

\bibitem{pmlr-v97-lee19e}
{\sc C.-P. Lee and S.~Wright}, {\em First-order algorithms converge faster than
  $o(1/k)$ on convex problems}, in Proceedings of the 36th International
  Conference on Machine Learning, vol.~97, 2019, pp.~3754--3762.

\bibitem{Lee16h}
{\sc J.~D. Lee, M.~Simchowitz, M.~I. Jordan, and B.~Recht}, {\em Gradient
  descent only converges to minimizers}, in 29th Annual Conference on Learning
  Theory, vol.~49, Columbia University, New York, New York, USA, 2016,
  pp.~1246--1257.

\bibitem{LessardRechtPackard2016}
{\sc L.~Lessard, B.~Recht, and A.~Packard}, {\em Analysis and design of
  optimization algorithms via integral quadratic constraints}, SIAM Journal on
  Optimization, 26 (2016), pp.~57--95,
  \url{https://doi.org/10.1137/15M1009597}.

\bibitem{LiuN89}
{\sc D.~C. Liu and J.~Nocedal}, {\em On the limited-memory {BFGS} method for
  large scale optimization}, Mathematical Programming, 45 (1989), pp.~503--528.

\bibitem{lu2025cupdlpx}
{\sc H.~Lu, Z.~Peng, and J.~Yang}, {\em {cuPDLPx}: A further enhanced
  {GPU}-based first-order solver for linear programming}, arXiv preprint
  arXiv:2507.14051,  (2025).

\bibitem{lu2024restarted}
{\sc H.~Lu and J.~Yang}, {\em Restarted {Halpern PDHG} for linear programming},
  arXiv preprint arXiv:2407.16144,  (2024).

\bibitem{Meh92a}
{\sc S.~Mehrotra}, {\em On the implementation of a primal-dual interior point
  method}, SIAM Journal on Optimization, 2 (1992), pp.~575--601.

\bibitem{mosek}
{\sc MOSEK ApS}, {\em The MOSEK Optimization Toolbox for MATLAB Manual},
  Copenhagen, Denmark, 2024, \url{https://docs.mosek.com}.
\newblock Version 10.2.

\bibitem{necoara2019linear}
{\sc I.~Necoara, Y.~Nesterov, and F.~Glineur}, {\em Linear convergence of first
  order methods for non-strongly convex optimization}, Mathematical
  programming, 175 (2019), pp.~69--107.

\bibitem{NemJLS09a}
{\sc A.~Nemirovski, A.~Juditsky, G.~Lan, and A.~Shapiro}, {\em Robust
  stochastic approximation approach to stochastic programming}, {SIAM} Journal
  on Optimization, 19 (2009), pp.~1574--1609.

\bibitem{NemY83}
{\sc A.~S. Nemirovski and D.~B. Yudin}, {\em Problem Complexity and Method
  Efficiency in Optimization}, John Wiley, 1983.

\bibitem{Nes83}
{\sc Y.~Nesterov}, {\em A method for unconstrained convex problem with the rate
  of convergence {$O(1/k^2)$}}, Doklady AN SSSR, 269 (1983), pp.~543--547.

\bibitem{Nes04}
{\sc Y.~Nesterov}, {\em Introductory Lectures on Convex Optimization: A Basic
  Course}, Springer Science and Business Media, New York, 2004.

\bibitem{nesterov2018lectures}
{\sc Y.~Nesterov}, {\em Lectures on Convex Optimization}, vol.~137 of Springer
  Optimization and its Applications, Springer, second~ed., 2018.

\bibitem{NesN93}
{\sc Y.~Nesterov and A.~S. Nemirovskii}, {\em Interior Point Polynomial Methods
  in Convex Programming}, SIAM Publications, Philadelphia, 1994.

\bibitem{NesP06a}
{\sc Y.~Nesterov and B.~T. Polyak}, {\em Cubic regularization of {Newton}
  method and its global performance}, Mathematical Programming, Series {A}, 108
  (2006), pp.~177--205.

\bibitem{Hogwild}
{\sc F.~Niu, B.~Recht, C.~{R\'e}, and S.~J. Wright}, {\em {\sc Hogwild!}: A
  lock-free approach to parallelizing stochastic gradient descent}, technical
  report, University of Wisconsin-Madison, June 2011.

\bibitem{NocW06}
{\sc J.~Nocedal and S.~J. Wright}, {\em {Numerical Optimization}}, Springer,
  New York, second~ed., 2006.

\bibitem{OneW18a}
{\sc M.~O'Neill and S.~J. Wright}, {\em Behavior of accelerated gradient
  methods near critical points of nonconvex problems}, Mathematical
  Programming, Series {B}, 176 (2019), pp.~403--427.

\bibitem{ParB13a}
{\sc N.~Parikh and S.~Boyd}, {\em Proximal algorithms}, Foundations and Trends
  in Optimization, 1 (2014), pp.~123--231.

\bibitem{PolR69}
{\sc E.~Polak and G.~{Ribi\`ere}}, {\em Note sur la convergence de {m\'ethodes}
  de directions conjug\'ees}, Revue {Fran\c{c}aise} d'Informatique et de
  Recherche {Op\'erationelle}, 16 (1969), pp.~35--43.

\bibitem{Pol64a}
{\sc B.~T. Polyak}, {\em Some methods of speeding up the convergence of
  iteration methods}, {USSR} Computational Mathematics and Mathematical
  Physics, 4.5 (1964), pp.~1--17.

\bibitem{Polyak1969}
{\sc B.~T. Polyak}, {\em {The conjugate gradient method in extremal problems}},
  USSR Computational Mathematics and Mathematical Physics, 9 (1969),
  pp.~94--112.

\bibitem{Pol87}
{\sc B.~T. Polyak}, {\em Introduction to Optimization}, Optimization Software,
  1987.

\bibitem{Powell76}
{\sc M.~J.~D. Powell}, {\em Some global convergence properties of a variable
  metric algorithm for minimization without exact line searches}, in Nonlinear
  Programming, SIAM-AMS Proceedings, Vol. IX, R.~W. Cottle and C.~E. Lemke,
  eds., SIAM Publications, 1976, pp.~53--72.

\bibitem{RobM51}
{\sc H.~Robbins and S.~Monro}, {\em A stochastic approximation method}, Annals
  of Mathematical Statistics, 22 (1951), pp.~400--407.

\bibitem{RoyOW18a}
{\sc C.~W. Royer, M.~O'Neill, and S.~J. Wright}, {\em A {Newton-CG} algorithm
  with complexity guarantees for smooth unconstrained optimization},
  Mathematical Programming, Series {A}, 180 (2020), pp.~451--488.

\bibitem{ShiDuJordanSu2022}
{\sc B.~Shi, S.~S. Du, M.~I. Jordan, and W.~J. Su}, {\em Understanding the
  acceleration phenomenon via high-resolution differential equations},
  Mathematical Programming,  (2022), pp.~1--70,
  \url{https://doi.org/10.1007/s10107-022-01864-w}.

\bibitem{Shub87}
{\sc M.~Shub}, {\em Global Stability of Dynamical Systems}, Springer, 1987.

\bibitem{smale2000algorithms}
{\sc S.~Smale}, {\em Algorithms for solving equations}, in The Collected Papers
  of Stephen Smale: Volume 3, World Scientific, 2000, pp.~1263--1286.

\bibitem{SpeT04}
{\sc D.~A. Spielman and S.-H. Teng}, {\em Smoothed analysis of algorithms: Why
  the simplex method usually takes polynomial time}, Journal of the Association
  for Computing Machinery, 51 (2004), pp.~385--463.

\bibitem{Ste83}
{\sc T.~Steihaug}, {\em The conjugate gradient method and trust regions in
  large scale optimization}, SIAM Journal on Numerical Analysis, 20 (1983),
  pp.~626--637.

\bibitem{Stu99}
{\sc J.~F. Sturm}, {\em Using {SeDuMi 1.02}, a {Matlab} toolbox for
  optimization over symmetric cones}, Optimization Methods and Software, 11
  (1999), pp.~625--653.

\bibitem{Sun21}
{\sc J.~Sun}, \url{https://sunju.org/research/nonconvex/}.

\bibitem{Tan87}
{\sc K.~Tanabe}, {\em Centered {N}ewton method for mathematical programming},
  in System Modeling and Optimization: Proceedings of the 13th {IFIP}
  conference, vol.~113 of Lecture Notes in Control and Information Systems,
  Berlin, August/September 1987 1988, Springer-Verlag, pp.~197--206.

\bibitem{TaylorHendrickxGlineur2017}
{\sc A.~B. Taylor, J.~M. Hendrickx, and F.~Glineur}, {\em Smooth strongly
  convex interpolation and exact worst-case performance of first-order
  methods}, Mathematical Programming, 161 (2017), pp.~307--345,
  \url{https://doi.org/10.1007/s10107-016-1009-x}.

\bibitem{Tod01a}
{\sc M.~J. Todd}, {\em Semidefinite optimization}, Acta Numerica, 10 (2001),
  pp.~515--560.

\bibitem{Todd2002}
{\sc M.~J. Todd}, {\em The many facets of linear programming}, Mathematical
  Programming, 91 (2002), pp.~417--436.

\bibitem{TodY90}
{\sc M.~J. Todd and Y.~Ye}, {\em A centered projective algorithm for linear
  programming}, Mathematics of Operations Research, 15 (1990), pp.~508--529.

\bibitem{sdpt3}
{\sc R.~H. {T\"ut\"unc\"u}, K.~C. Toh, and M.~J. Todd}, {\em {SDPT3}---A
  {Matlab} software package for semidefinite-quadratic-linear programming,
  version {3.0}}, August 2001.

\bibitem{TutTT03}
{\sc R.~H. {T\"ut\"unc\"u}, K.~C. Toh, and M.~J. Todd}, {\em Solving
  semidefinite-quadratic-linear programs using {SDPT3}}, Mathematical
  Programming, Series {B}, 95 (2003), pp.~189--217.

\bibitem{van2020deterministic}
{\sc J.~van~den Brand}, {\em A deterministic linear program solver in current
  matrix multiplication time}, in Proceedings of the Thirty-First Annual
  ACM-SIAM Symposium on Discrete Algorithms, SODA '20, USA, 2020, Society for
  Industrial and Applied Mathematics, p.~259–278.

\bibitem{Vershynin2009}
{\sc R.~Vershynin}, {\em Beyond hirsch conjecture: walks on random polytopes
  and smoothed complexity of the simplex method}, SIAM Journal on Computing, 39
  (2009), pp.~646--678, \url{https://doi.org/10.1137/070683386}.

\bibitem{weiaccelerated25}
{\sc J.~Wei and L.~Chen}, {\em Accelerated over-relaxation heavy-ball method:
  Achieving global accelerated convergence with broad generalization}, in The
  Thirteenth International Conference on Learning Representations, 2025.

\bibitem{Wri97}
{\sc S.~J. Wright}, {\em Primal-Dual Interior-Point Methods}, SIAM,
  Philadelphia, PA, 1997.

\bibitem{WriNF08a}
{\sc S.~J. Wright, R.~D. Nowak, and M.~A.~T. Figueiredo}, {\em Sparse
  reconstruction by separable approximation}, IEEE Transactions on Signal
  Processing, 57 (2009), pp.~2479--2493.

\bibitem{WrightRecht22}
{\sc S.~J. Wright and B.~Recht}, {\em Optimization for Data Analysis},
  Cambridge University Press, May 2022.

\bibitem{xiong2025high}
{\sc Z.~Xiong}, {\em High-probability polynomial-time complexity of restarted
  {PDHG} for linear programming}, arXiv preprint arXiv:2501.00728,  (2025).

\bibitem{young1953richardson}
{\sc D.~Young}, {\em On {Richardson}'s method for solving linear systems with
  positive definite matrices}, Journal of Mathematics and Physics, 32 (1953),
  pp.~243--255.

\bibitem{YudinNemirovski1976}
{\sc D.~B. Yudin and A.~S. Nemirovski}, {\em Informational complexity and
  efficient methods for the solution of convex extremal problems},
  {\`E}konomika i Matematicheskie Metody, 12 (1976), pp.~357--369.
\newblock English translation in Matekon, vol. 13(2), 1980, pp. 3--25.

\end{thebibliography}
\end{document}